%% file: RenewalScheduling.tex
\documentclass[a4paper,british,english]{article}
\usepackage[T1]{fontenc}
\usepackage[utf8]{inputenc}
\usepackage{xcolor}
\usepackage{array}
\usepackage{mathtools}
\usepackage{amsmath}
\usepackage{amssymb}
\usepackage{graphicx}

\makeatletter

\pdfpageheight\paperheight
\pdfpagewidth\paperwidth

\providecommand{\tabularnewline}{\\}

\usepackage{tikz}
\usepackage{hyperref}

\makeatother

\usepackage{babel}
\usepackage[style=authoryear]{biblatex}
\begin{document}
\title{An optimization model for renewal scheduling and traffic flow routing}
\author{Tomas Lidén, Martin Aronsson}
\maketitle
\begin{abstract}
This document is the third sub-report from the research project SATT
(Samplanering av trafikpåverkande åtgärder och trafikflöden, modellstudie
/ Coordinated planning of temporary capacity restrictions and traffic
flows, model study). The report summarizes the modelling alternatives
that have been considered. The main result of the model study is that
we propose a bi-level approach for handling the two parts of the problem,
i.e. the scheduling of renewal projects and the subsequent adjustments
of railway traffic flow. The contributions are: (1) a mixed integer
formulation for scheduling project tasks, and (2) a network flow formulation
for establishing the possible volumes of railway traffic during the
imposed capacity restrictions of these project tasks. The intended
usage of the model is for economical planning of a specific production
year, which takes place before the timetabling process starts. This
document introduces the problem and gives a full description of (1)
along with results from computational experiments that have been conducted
with this model. An accompanying report gives the details of the final
formulation of (2).
\end{abstract}

\section{Introduction}

An infrastructure manager is responsible for the capacity planning
of a transportation network. The long term planning consists of deciding
on renewal, upgrade, and major maintenance activities which constitute
large budget volumes for a national railway system. Thus, a multi-annual
economical planning is necessary. This planning is done before the
details of the railway traffic (aka the timetable) are known. Nevertheless,
there is a need to understand the implications that the project planning
has on the train traffic that will take place alongside the project
tasks. Hence, the long-term economical planning should both ensure
that the considered infrastructure projects can be completed, while
at the same time limiting the negative impact of these projects on
the on-going traffic. This planning problem is sometimes labeled as
``intervention planning'' \parencite{burkhalter_novel_2021}.

In this work we study the coordination of major capacity affecting
work activities and traffic flows, by developing a method for optimal
scheduling of project tasks, which considers both the project aspects
and the traffic throughput. We are particularly interested in major
work tasks that take place on an existing railway network, such as
upgrades, renewals and larger maintenance activities\footnote{From here on, we make no particular distinction between the terms
``renewal'' and ``upgrade'' projects. Note also that major maintenance
activities are included in the project tasks that we consider.}. Thus we disregard investment projects and the construction of new
infrastructure which does not interact with any on-going train traffic.
Also, we do not consider activities with minor capacity impact, such
as shorter and recurrent regular maintenance, which are, or can be,
handled in the tactical and operational planning.

\subsection{Background, motivation and problem description\label{subsec:Background-problem-description}}

The current planning of railway renewal projects in Sweden is done
with a limited view of the traffic impact. Each project is treated
individually and network wide effects, particularly the combinatorial
effects of several or all projects, are not supported by the available
planning methods or tools. In fact, very few examples of such planning
support can be found in the research literature. Instead, rough assumptions
are made regarding traffic volumes, and the capacity impact is estimated
based on experience and with advisory tables regarding which work
types and work locations that might or should not be conducted at
the same time. There are few performance indicators and metrics for
measuring the traffic impact, limited knowledge about the actual capacity
effects of previously conducted projects, and no strong decision support
for dimensioning, configuring and scheduling renewal projects.

On the other hand, there is a growing awareness for the need of well
founded quantitative planning methods. Especially when governmental
spending is increasing as a consequence of aging infrastructure, historical
asset management deficits, and an increase in railway transportation
demand. In this situation, the coordination of capacity needs from
renewal and traffic activities must be balanced with the intention
of achieving an optimal outcome for society (according to available
resources)\footnote{Swedish legislation currently stipulates that all infrastructure planning
should be based on the principle of socioeconomic efficiency, while
European directives strive to maximize societal benefits.}.

The need for improved decision support regarding the scheduling of
renewal projects in the long-term economical (intervention) planning
is the motivation for this study. The purpose of the work is to investigate
the possibility of using mathematical optimization techniques for
this planning stage.

The requirement details for the planning problem has been elaborated
in \textcite{liden_krav_2021}. The resulting problem definition can
be stated as follows: Given a railway network with an assumed transportation
flow demand and a set of considered renewal projects, each consisting
of a sequence of work tasks with specific resource needs (e.g. machines
and crew), schedule as many projects as possible with a minimal combined
cost for task resource usage and traffic disruptions, such that limitations
regarding project timing, resource availability, capacity restrictions,
and traffic flow, is respected. The intention is to be able to solve
problem instances for national networks and a scheduling period of
one operational year.

\subsection{Literature overview}

Scheduling methods for railways have traditionally focused on train
traffic. Some examples of well studied problems are timetabling, rolling
stock and crew scheduling, rostering, routing, rescheduling, speed
optimization, etc. In recent years the coordination and joint scheduling
of maintenance and traffic have attracted some research activity,
but so far mostly regarding tactical planning where each train service
is individually scheduled. Such optimization model have been able
to solve weekly problem instances \parencite{liden_utformning_2020}.
Although the basic problem type (coordination of work tasks and train
traffic) is similar, these models are too fine-grained regarding the
traffic representation to suit the requirements for yearly renewal
scheduling. However, the aggregated network representation and capacity
control applied for example by \textcite{liden_optimization_2017}
is an approach that is suitable in our problem setting.

As for aggregated handling of national railway traffic flows, there
are limited references in the optimization literature. However, in
\textcite{aronsson_ttjob_2019} a model for describing railway transportation
in the form of so called ``transport service classes'' is presented.
The idea is to describe a (future) national transportation plan with
these transport service classes, based on origin-destination relations,
travel time and crucial service requirements, while giving the freedom
of fulfilling these transport service classes with different routing
options. When traffic is described in this form, there is no need
to consider individual trains or train paths. Instead the traffic
can be represented as flows, where appropriate link capacity constraints
should ensure that a feasible timetable exists for this flow. This
kind of traffic representation is suitable for our problem setting.

The above given references was the starting inspiration for this work.
In \textcite{liden_krav_2021} a survey of other related research
is presented. The literature search focused on publications that treat
both track works and train traffic in a long-term setting, where the
latter is represented on an aggregated level---not as timetables
and individual trains but as flows of traffic. The publications were
categorized according to the temporal resolution regarding work and
traffic flow. Models where the activity type are shorter than the
time period lengths, such that they always can start and complete
within one time period, we label as ``static''\footnote{Note that the term ``static'' only says that each activity will
affect one and only one time period. The actual scheduling values
(number of work tasks, traffic flow, etc) are still variable within
each time period.}. If the activities span over more than one time period, such that
the partial overlapping of activities becomes crucial, we use the
term ``dynamic'' (see for example Skutella (2009) regarding timed,
or dynamic, network flows). Thus, we get four types of models, namely:
\begin{itemize}
\item Coarse models, where both work activities and traffic flow is static.
A common example is to consider which year to perform different projects,
but where the project time is considerably shorter and where traffic
impact is treated without consideration of synchronization or network
effects.
\item Work-flexible models with static traffic, where track work activities
span over several time periods but the traffic flows have no interaction
between time periods. A typical time resolution for such models would
be to use a weekly or daily period length.
\item Traffic-flexible models with static work activities, where traffic
flows will interact between time periods while each work activity
only affect one time period. Thus there is no partial overlap between
work activities and the capacity restrictions imposed on the traffic
is fixed within each time period.
\item Work- and traffic-flexible models, where both track work activities
and traffic flows span over more than one time period, such that the
dynamic effects over the network topology must be considered.
\end{itemize}
We aim at the last type of model, but the bulk of research literature
has only been found in the two first types (with static traffic flow
modeling). From the coarse models we have not found any applicable
modeling inspiration apart from \textcite{li_biobjective_2020} which
use different time scales for track works and link status. This is
an idea that could be applied to our case for upgrade projects and
traffic flow.

The work-flexible models with static traffic show a large variation
regarding level of details and which aspects that are considered for
the track works. However, all these models have the limitation that
each track work only consists of one uninterrupted activity (or task).
These references however indicate that a combination of a mixed integer
programming and a network flow formulation (for track works and traffic
respectively) might be suitable in our problem setting.

Finally we have not been able to find any traffic-flexible models
based on traffic flows with properties that fit railway systems and
the planning problem we consider. This indicates either that our search
has been too limited or that we have found a research gap where important
and interesting contributions can be made.

\subsection{Outline of work and document overview}

The modeling work started by considering an approach based on a set
partitioning formulation with two types of columns, one for project
schedules and one for traffic flow routing. A model sketch of this
is given in Appendix \ref{app:Set-partition-model}, first for an
integrated model, then for a bi-level approach. At this stage focus
was given to the traffic model (the lower level of the bi-level approach)
and instead of column generation we worked on network flow formulations.
Three variants (with different layer definitions) were studied. An
overview of these alternatives is given in Section \ref{subsec:Flow-model-overview}
and the details of the first two can be found in Appendix \ref{app:Flow-models}.
The third formulation of the traffic flow model is documented in the
accompanying report \textcite{aronsson_flows_2021}. The last stage
of the modeling work concerned the scheduling of upgrade projects
(the upper level of the bi-level approach), where a mixed integer
formulation was developed. The mathematical model is presented in
Section \ref{subsec:MIP-upper-level} followed by results from computational
experiments in Chapter \ref{sec:Computational-experiments}. The document
ends with some final conclusions and discussion in Chapter \ref{sec:Conclusions}.

\section{Mathematical formulation}

\begin{figure}
\begin{centering}
\input{figures/bilevel_tikz.tex}
\par\end{centering}
\caption{\label{fig:bi-level}Bi-level model outline}
\end{figure}
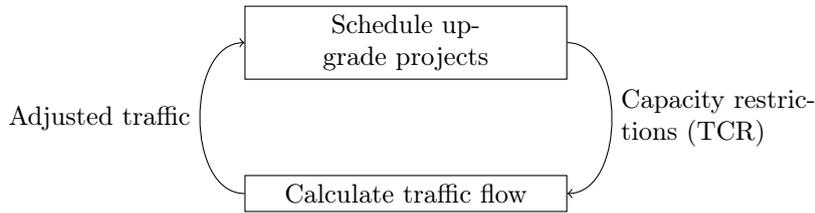

We propose a bi-level approach where the two parts of the problem
are treated in two interconnected models, as illustrated in Figure
\ref{fig:bi-level}. The upper level is the scheduling of the upgrade
projects according to an initial cost function. For each operational
day which is affected by (temporary) capacity restrictions (TCRs),
the subsequent traffic adjustments are calculated by the lower level
problem. Thus, only the affected traffic days need to be considered
in the lower level, which reduces the problem size. The obtained traffic
adjustments can then be considered by the upper level problem in order
to find a better TCR schedule, taking both the cost for the upgrade
projects and the traffic routing into account. The two processes are
iterated until some convergence or optimality criteria is met.

In the remainder of this chapter we describe the modeling in more
detail. Section \ref{subsec:Basic-notation} gives some basic notation,
followed by an outline of the flow-based formulation models that have
been investigated for the lower level problem (Section \ref{subsec:Flow-model-overview}).
Finally, the complete formulation of the upper level model is given
in Section \ref{subsec:MIP-upper-level}.

\subsection{Basic notation\label{subsec:Basic-notation}}

The following basic notation will be used throughout the document.

\subsubsection*{Planning period}

\begin{tabular}{cl}
$T$ & a sequence of time periods $t$\tabularnewline
$\beta_{t},\delta_{t}$ & start and length of time period $t\in T$\tabularnewline
\end{tabular}

\noindent As a simplification we will mostly assume a unit time period
length ($\delta_{t}=1$) where each time period $t$ starts at $\beta_{t}=t-1$,
and spans the open-ended time interval $[t-1,t)$.

\subsubsection*{Railway network}

\begin{tabular}{cl}
$N$ & a set of network nodes $n$\tabularnewline
$L$ & a set of network links $l$\tabularnewline
$C_{l}^{nom}$ & the nominal capacity (number of trains per time unit) for link $l\in L$\tabularnewline
\end{tabular}

\subsubsection*{Wanted traffic}

\begin{tabular}{c>{\raggedright}p{0.9\columnwidth}}
$H$ & a set of train types $h$\tabularnewline
$R$ & a set of traffic relations $r$, each defined by an origin-destination-train
type tuple $(o_{r},d_{r},h_{r})\in N\times N\times H$\tabularnewline
\end{tabular}

\noindent For each traffic relation there is a limited set of routing
options through the network, which we will discuss in subsequent sections.

\subsubsection*{Renewals / upgrade projects}

\begin{tabular}{cl}
$U$ & a set of upgrade projects $u$, which in turn consists of a sequence
of tasks.\tabularnewline
\end{tabular}

\noindent The details regarding project tasks will be further described
in Section \ref{subsec:MIP-upper-level}.

\subsection{Flow based approach of the lower level problem\label{subsec:Flow-model-overview}}

The idea of the flow based formulation is to describe the problem
as a layered space-time network problem, where time is discretized
into relatively coarse time periods, say of length 1--6 hours, such
that the link travel time for all trains are shorter than any time
period length. If we let $\nu$ be the travel time expressed as a
share of the time period length, this means that we require $0<\nu<1$.
Assuming that flow is evenly distributed over the time period, then
$1-\nu$ of the flow will arrive within the same time period and $\nu$
of the flow in the next time period. We further demand flow continuity
such that no incoming flow from the previous time period ($t-1$)
into time period $t$ is allowed to proceed to the succeeding time
period ($t+1$). With this requirement we avoid situations where flow
might ``jump hurdles'', and in addition the formulation complexity
is reduced. For a situation of two connected links $p$ and $l$,
with (relative) travel times $\nu_{p}$ and $\nu_{l}$, the requirement
is then $\nu_{p}+\nu_{l}\leq1$ (or $1-\nu_{p}-\nu_{l}\geq0$). A
simple and conservative limitation is to require $\nu\leq\frac{1}{2}$,
but a more flexible approach is to surround (long) links where $\nu\geq0.5$
with shorter ones, such that the requirement holds.

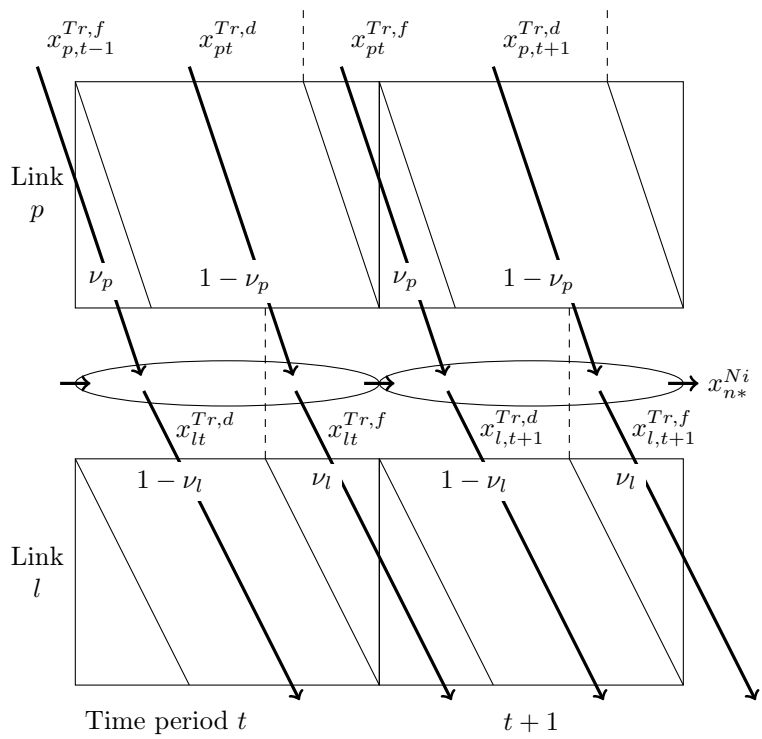
\begin{figure}
\begin{centering}
\input{figures/flowidea_tikz.tex}
\par\end{centering}
\caption{\label{fig:flow-idea}Illustration of flow relations.}
\end{figure}

These concepts are illustrated in Figure \ref{fig:flow-idea} for
two successive links ($p$ and $l$), where bold arrows indicate flow
variables, named as follows:

\noindent %
\begin{tabular}{c>{\raggedright}p{0.9\columnwidth}}
\noalign{\vskip3pt}
$x_{lt}^{Tr,d}$ & Transportation flow over link $l$, starting in time period $t$ and
arriving in the same time period (direct flow)\tabularnewline
\noalign{\vskip3pt}
$x_{lt}^{Tr,f}$ & Transportation flow over link $l$, starting in time period $t$ and
arriving in the next period $t+1$ (forward flow)\tabularnewline
\noalign{\vskip3pt}
$x_{nt}^{Ni}$ & Node inventory flow for node $n$ between time period $t$ and $t+1$\tabularnewline[3pt]
\noalign{\vskip3pt}
\end{tabular}

The consequence of the continuity requirement\footnote{when considering transportation flow that can traverse as fast as
possible, i.e. there is no inventory flow or any link capacity restrictions} is that $x_{lt}^{Tr,f}\leq x_{pt}^{Tr,d}$, and that we can impose
constraints that control the volumes of direct and forward flow, so
as to achieve proper travel times through the network.

The next important question concerns the layer definition. We have
investigated three approaches:
\begin{itemize}
\item One layer per train type $h$
\item One layer per origin-destination relation $r$
\item One layer per routing option $p$ of each traffic relation $r$, where
each routing option constitutes one path through the network
\end{itemize}
The two first options are described in more detail in Appendix \ref{app:Flow-models}
along with some computational experiments. Both models allow for flow
to split and merge within each layer. The first model does not properly
handle reachability requirements for the traffic relations, while
the second (which is a multi-commodity network flow model with side
constraints) resolves this issue. However, both models have a major
drawback when flow runs over many links, which results in small fractional
flow volumes being pushed forward in time (due to the side constraints
that control the volumes of direct and forward flow). We have not
found any good solution to this problem, and have therefore developed
the third model, with one layer for each relation+routing option.
This approach gives a planar graph for each layer, and enables another
way of imposing travel time restrictions by giving $\nu_{rpn}$ values
that state the shortest (relative) travel time (from the origin) for
reaching node $n$ of path $p$ and relation $r$, which in turn will
restrict the maximum flow volumes that can reach a node at each time
period $t$. This model correctly pushes traffic forward in time without
the drawbacks of the first two models, although the number of layers,
variables and constraints will increase.

The full description of the ``routing-option flow model'' is given
in the accompanying report \textcite{aronsson_flows_2021}.

\subsection{MIP formulation of the upper level problem\label{subsec:MIP-upper-level}}

This section contains a model write-up of the upper level problem
as a MIP formulation. First, we introduce some additional input data:

\noindent %
\begin{tabular}{l>{\raggedright}p{0.8\columnwidth}}
\noalign{\vskip3pt}
$T_{u}\subseteq T$ & The allowed set of time periods (from start to end) for upgrade project
$u$\tabularnewline
\noalign{\vskip3pt}
$W_{u}$ & The number of work tasks in upgrade project $u$\tabularnewline
\noalign{\vskip3pt}
$\eta_{uw}$ & The length (number of time periods) for work task $w$ of upgrade
project $u$\tabularnewline
\noalign{\vskip3pt}
$\gamma_{uw},\theta_{uw}$ & The min and max length (number of time periods) of rest time between
work task $w$ and $w+1$ of upgrade project $u$\tabularnewline
\noalign{\vskip3pt}
$b_{uwl}$ & The number of capacity blockings caused by work task $w$ of upgrade
project $u$ on link $l$\tabularnewline
\noalign{\vskip3pt}
$F_{d}(\mathbf{x}^{d})$ & Disruption cost for traffic assignment $\mathbf{x}^{d}$ of day $d$\tabularnewline
\noalign{\vskip3pt}
$K$ & The set of all work resources $k$\tabularnewline
\noalign{\vskip3pt}
$g_{uwk}$ & The amount of resource usage for work resource $k$ when performing
work task $w$ of project $u$\tabularnewline
\noalign{\vskip3pt}
$Q_{k}$ & The available quantity of work resource $k$\tabularnewline[3pt]
\noalign{\vskip3pt}
\end{tabular}

\noindent The basic decision variables\footnote{\noindent The binary restriction on either the $y$ or the $z$ variables
can be omitted, as the constraints will enforce integrality of the
other.} in the model are

\noindent %
\begin{tabular}{l>{\raggedright}p{0.7\columnwidth}}
\noalign{\vskip3pt}
$y_{uwt}\in\mathbb{Z}_{2}^{U\times W\times T}$ & (Binary) Whether work task $w$ of upgrade project $u$ is taking
place during time period $t$ or not\tabularnewline
\noalign{\vskip3pt}
$z_{uwt}\in\mathbb{Z}_{2}^{U\times W\times T}$ & (Binary) Whether work task $w$ of upgrade project $u$ is starting
in time period $t$ or not\tabularnewline
\noalign{\vskip3pt}
$s_{u}\in\mathbb{Z}_{2}^{U}$ & (Binary) Whether upgrade project $u$ is canceled (not scheduled)
or not\tabularnewline
\noalign{\vskip3pt}
$q_{k}\in\mathbb{R}_{+}^{K}$ & The max usage of resource $k$ in any time period\tabularnewline[3pt]
\noalign{\vskip3pt}
\end{tabular}

\noindent We assume that each work task has a fixed duration (can
easily be changed) and that the tasks within an upgrade project should
be performed in their given order. The decision freedom lies in selecting
start times of each work task, while respecting the spacing requirements
between successive tasks. The optimization problem can then be formulated
as follows:

\begin{multline*}
\textrm{minimize}\sum_{u\in U}c_{u}^{U}s_{u}+\sum_{k\in K}c_{k}^{K}q_{k}+\sum_{d\in WD\left(\bar{\mathbf{y}}\right)}F_{d}\left(\mathbf{x}^{d}\right)\\
+\sum_{l\in L.t\in T}\lambda_{lt}\left(a_{lt}\left(\mathbf{x}\right)+\sum_{u\in U}\sum_{w=1}^{W_{u}}b_{uwl}y_{uwt}-C_{l}^{nom}\right)
\end{multline*}

\noindent subject to constraints for:

scheduling each work task once unless project is canceled

\begin{align}
\forall u\in U,w\in\left[1,W_{u}\right]: & \sum_{t\in T_{u}}z_{uwt}=1-s_{u};\label{eq:schedule-every-task}
\end{align}

linking on-going work variables to work start variables

\begin{align}
\forall u\in U,w\in\left[1,W_{u}\right],t\in T_{u}: & z_{uwt}\geq y_{uwt}-y_{uw,t-1};\label{eq:link-y-and-z}
\end{align}

length of each work task

\begin{align}
\forall u\in U,w\in\left[1,W_{u}\right],t\in T_{u}: & \sum_{i=1}^{\min(\eta_{uw},t)}z_{uw,t-i+1}=y_{uwt};\label{eq:task-length}
\end{align}

no overlapping tasks within a project

\begin{align}
\forall u\in U,t\in T_{u}: & \sum_{w=1}^{W_{u}}y_{uwt}\leq1;\label{eq:no-task-overlap}
\end{align}

successive order of tasks within a project

\begin{align}
\forall u\in U,w\in\left[1,W_{u}-1\right],i\in\left[1,W_{u}-w\right],t\in T: & \sum_{j=1}^{t}\left(z_{uwj}-z_{u,w+i,j}\right)\geq0;\label{eq:task-order}
\end{align}

minimum and maximum length of rest time between tasks

\begin{align}
\forall u\in U,w\in\left[1,W_{u}-1\right],t\in T_{u}: & \sum_{i=1}^{\min(\gamma_{uw},\left|T\right|-t)}z_{u,w+1,t+i}\leq\sum_{j=1}^{t}z_{uwj}-y_{uwt};\label{eq:min-rest}
\end{align}

\begin{align}
\forall u\in U,w\in\left[1,W_{u}-1\right],t\in T_{u}: & \sum_{i=1}^{\min(\theta_{uw},\left|T\right|-t)}z_{u,w+1,t+i}\geq\sum_{j=1}^{t}\left(z_{uwj}-z_{u,w+1,j}\right)-y_{uwt};\label{eq:max-rest}
\end{align}

and limited resource usage

\begin{align}
\forall k\in K,t\in T: & \sum_{u\in U}\sum_{w=1}^{W_{u}}g_{uwk}y_{uwt}\leq q_{k};\label{eq:resource-use}
\end{align}

\begin{align}
\forall k\in K: & q_{k}\leq Q_{k}.\label{eq:resouce-limit}
\end{align}

Note that the objective function penalizes unscheduled work tasks,
resource usage, traffic disruption cost and a Lagrangian relaxation
term for line capacity utilization. Work costs could be added, and
coordination of certain activities (e.g. interlocking or traffic management
system loading) can be achieved by setting time period limitations
on work task level (i.e. using data sets $T_{uw}$ rather than $T_{u}$).

This formulation uses a large number of variables and constraints.
As an example, consider a one year planning period divided into 6
hour periods, which results in $365\times4=1460$ time periods. Then
the scheduling of 100 upgrade projects, each consisting of (in average)
10 work tasks, where each project has an allowed scheduling period
of 4 months ($\Rightarrow4\times30\times4=480$ time periods) will
require in the order of 1 M variables and 2.5 M constraints. However,
if time period limitations are set on work task level, such that each
work task has on average a one month scheduling window, the problem
size is reduced to about 250 k variables and 750 k constraints. Such
problem sizes can be tractable for modern LP and MIP solvers. An alternative
could be to use a CP approach, especially for handling resource limitations
and for further pruning of the variable domains (while still handling
the optimization problem in a LP/MIP solver). It could also be considered
to use real/integer valued task start variables (rather than the binary
$z_{uwt}$) so as to reduce the number of variables and get a more
compact formulation.

The objective function relies on knowledge of the traffic disruption
cost ($F_{d}$) and the traffic capacity usage ($a_{lt}$) over each
link and time period, as calculated by the lower level traffic assignment
problem. Without knowledge of the actual traffic assignment (for a
particular project plan), it is necessary to have some estimate of
the traffic impact for a certain project schedule. We can use the
blocking values $b_{uwl}$ to setup a proxy for the traffic impact.
As a lower bound we can use the largest blocking value for any link
during each time period---assuming perfect coordination such that
all other blocking values during the same time period does not further
affect the traffic. As an upper bound we instead sum all blocking
values used at all links and time periods, which corresponds to having
no traffic coordination at all. To establish these values we introduce
the variables

\noindent %
\begin{tabular}{l>{\raggedright}p{0.6\columnwidth}}
\noalign{\vskip3pt}
$x_{t}^{b,m}\in\mathbb{R}_{+}^{T}$ & The largest blocking value for any project on any link in time period
$t$\tabularnewline
\noalign{\vskip3pt}
$B^{LB},B^{UB}\in\mathbb{R}_{+}$ & The lower and upper bound values for number of blocked (affected)
trains\tabularnewline[3pt]
\noalign{\vskip3pt}
\end{tabular}

\noindent and the constraints

\begin{align}
\forall u\in U,l\in L,t\in T: & \sum_{w=1}^{W_{u}}b_{uwl}y_{uwt}\leq x_{t}^{b,m}\label{eq:blocking-count}
\end{align}

\begin{equation}
B^{LB}=\sum_{t\in T}x_{t}^{b,m}\label{eq:min-blocking}
\end{equation}

\begin{equation}
B^{UB}=\sum_{u\in U}\sum_{w=1}^{W_{u}}\sum_{l\in L}\sum_{t\in T}b_{uwl}y_{uwt}\label{eq:max-blocking}
\end{equation}

Note that $x_{t}^{b,m}$ will become the largest blocking value of
any upgrade project $u$ that affect the same link in the same time
period. Thus it is assumed that project tasks on the same link are
possible to coordinate. (A more conservative approach would be to
sum up overlapping blocking values from different projects in the
left hand side of the first constraint.) This model can be further
improved by including a link-time-dependent cost factor, such that
blockings can be made more or less costly depending on the normal
traffic load over specific links and time periods. Thus, the values
$x_{t}^{b,m}$ and $B^{LB/UB}$ will measure the estimated traffic
disruption rather than unavailable train slots.

The values $B^{LB}$ and $B^{UB}$ can be used in the objective function,
for example by using a parameter value, $0\leq\zeta\leq1$ that estimate
the level of possible traffic coordination, and adding a cost component

\[
c^{disr}\left[\zeta B^{LB}+\left(1-\zeta\right)B^{UB}\right]
\]

Additionally, we might want to track time periods with disturbed traffic,
using a variable $x_{t}^{d}\in\mathbb{Z}_{2}^{T}$ and the following
constraints

\begin{align*}
\forall u\in U,w\in\left[1,W_{u}\right],t\in T: & y_{uwt}\leq x_{t}^{d}
\end{align*}

\begin{align*}
\forall t\in T: & x_{t}^{d}\leq\sum_{u\in U}\sum_{w=1}^{W_{u}}y_{uwt}
\end{align*}

\noindent with the possibility to include $x_{t}^{d}$ in the objective
function.

For the separation of work tasks it is possible to include an additional
set of variables

\noindent %
\begin{tabular}{l>{\raggedright}p{0.6\columnwidth}}
\noalign{\vskip5pt}
$e_{uw}\in\mathbb{R}_{+}^{U\times W}$ & Start time for work task $w$ of upgrade project $u$\tabularnewline[5pt]
\end{tabular}

\noindent and the following constraints

\begin{align}
\forall u\in U,w\in\left[1,W_{u}\right]: & e_{uw}=\sum_{t\in T}tz_{uwt}\label{eq:event-time}\\
\forall u\in U,w\in\left[1,W_{u}-1\right]: & e_{u,w+1}-e_{uw}\geq\left(1-s_{u}\right)\left(\eta_{uw}+\gamma_{uw}\right)\label{eq:min-separation}\\
\forall u\in U,w\in\left[1,W_{u}-1\right]: & e_{u,w+1}-e_{uw}\leq\left(1-s_{u}\right)\left(\eta_{uw}+\theta_{uw}\right)\label{eq:max-separation}
\end{align}

Computational experiments on a limited set of instances indicate that
these variables and constraints improve solving performance in most
cases. Furthermore these constraints make it possible to disable one
or both of the previously given constraints (\ref{eq:min-rest}) and
(\ref{eq:max-rest}) for minimum and maximum length of rest time between
tasks. Our experiments indicate a slight performance gain for removing
the minimum rest time constraints (\ref{eq:min-rest}).

\section{Computational experiments\label{sec:Computational-experiments}}

The MIP formulation for the upper level problem (upgrade project scheduling)
has been implemented in MiniZinc (version 2.5.5) and solved with COIN-BC
(version 2.10.5) as optimization solver. The experiments have been
run on a HP EliteBook 830 G5 with an 1.6 GHz i5 processor, 8 GB RAM,
and 64 bit Windows 10 Enterprise OS (version 19043.867). We include
the event time variables $e_{uw}$ and the corresponding task separation
constraints (\ref{eq:event-time})--(\ref{eq:max-separation}), but
exclude the minimum rest time constraints (\ref{eq:min-rest}) and
(\ref{eq:max-rest}). A small data instance have been used, with 8
nodes and 9 links (see Figure \ref{fig:Small-flow-case}). Three projects
with 3, 7 and 5 work tasks respectively are to be scheduled and the
task lengths, min/max rest time, resource usage and link blocking
values are given in Table \ref{tab:Project-settings}. Directly from
the input data we can deduce that the projects must span 14--28,
34--58 and 28--55 time periods respectively, that the minimum resource
usage is 1 and 10 (for resources R1 and R2), and that the minimum
link blocking for each project is 20, 38 and 52 respectively.

\begin{table}
\begin{centering}
\begin{tabular}{|c|c|c|cc|ccccccccc|}
\hline 
Proj/ & Task & Rest time & \multicolumn{2}{c|}{Usage} & \multicolumn{9}{c|}{Link blocking}\tabularnewline
Task & len. & Min/Max & R1 & R2 & A-C & B-C & C-E & C-F & E-H & F-H & D-E & E-F & F-G\tabularnewline
\hline 
\hline 
P1/W1 & 2 & 2/6 & 0 & 4 &  &  & 2 &  & 2 &  &  &  & \tabularnewline
\textcolor{lightgray}{P1}/W2 & 4 & 4/8 & 1 & 8 & 1 &  & 3 &  & 3 &  &  &  & \tabularnewline
\textcolor{lightgray}{P1}/W3 & 2 & -/- & 0 & 4 &  &  & 2 &  & 2 &  &  &  & \tabularnewline
\hline 
P2/W1 & 2 & 2/6 & 0 & 4 & 2 &  &  &  &  &  &  &  & \tabularnewline
\textcolor{lightgray}{P2}/W2 & 2 & 2/6 & 0 & 4 &  & 2 &  &  &  &  &  &  & \tabularnewline
\textcolor{lightgray}{P2}/W3 & 2 & 2/6 & 0 & 6 & 1 & 1 & 1 & 1 &  &  &  &  & \tabularnewline
\textcolor{lightgray}{P2}/W4 & 4 & 4/8 & 1 & 10 & 2 &  & 3 &  &  &  &  &  & \tabularnewline
\textcolor{lightgray}{P2}/W5 & 4 & 4/8 & 1 & 10 &  & 2 &  & 3 &  &  &  &  & \tabularnewline
\textcolor{lightgray}{P2}/W6 & 2 & 2/6 & 0 & 6 & 1 & 1 & 1 & 1 &  &  &  &  & \tabularnewline
\textcolor{lightgray}{P2}/W7 & 2 & -/- & 0 & 6 & 1 & 1 & 1 & 1 &  &  &  &  & \tabularnewline
\hline 
P3/W1 & 1 & 2/7 & 0 & 4 &  &  &  &  &  &  & 2 &  & 2\tabularnewline
\textcolor{lightgray}{P3}/W2 & 2 & 2/8 & 0 & 8 &  &  &  &  &  &  & 1 & 2 & 1\tabularnewline
\textcolor{lightgray}{P3}/W3 & 8 & 8/20 & 1 & 8 &  &  &  &  &  &  & 5 & 5 & 5\tabularnewline
\textcolor{lightgray}{P3}/W4 & 2 & 2/6 & 0 & 8 &  &  &  &  &  &  & 1 & 2 & 1\tabularnewline
\textcolor{lightgray}{P3}/W5 & 1 & -/- & 0 & 4 &  &  &  &  &  &  & 2 &  & 2\tabularnewline
\hline 
\end{tabular}
\par\end{centering}
\caption{\label{tab:Project-settings}Project schedule settings for small test
case}
\end{table}

In the experiments we will use unit resource costs ($c_{k}^{K}=1,\forall k)$,
the traffic disruption cost factor $c^{disr}=0.1$ and a traffic coordination
parameter value $\zeta=0.75$. The solving performance will depend
on the scheduling freedom for each project, but we will not impose
any restrictions on the time windows for each project. Thus there
is a large number of symmetric solutions.

Three sets of experiments will be reported. In Section \ref{subsec:Traffic-disruption-cost}
and Section \ref{subsec:Resource-usage} we will focus solely on the
traffic disruption cost and the resource usage respectively. From
these experiments we obtain lower bounds which are then used in Section
\ref{subsec:Combined-traffic-and-resources} to \foreignlanguage{british}{analyse}
the combined cost of both traffic disruption and resource cost.

\subsection{Traffic disruption cost\label{subsec:Traffic-disruption-cost}}

In this experiment we find lower bounds for the traffic disruption
cost by forcing all projects to be scheduled ($s_{u}=0,\forall u$),
disregarding the resource costs ($c_{k}^{K}=0,\forall k$), and solving
the scheduling problem for four planning period sizes ($\left|T\right|=35,40,45,50$).
The results are given in Table \ref{tab:Traffic-disruption-lb}, which
reports the resulting problem sizes (number of variables x constraints),
time for finding the integer solution (seconds), total solution time
(seconds), the lower bound blocking value ($B^{LB}$), the disruption
cost, and the resource usage (R1+R2) for the respective schedule solutions.
The latter values are given in parentheses to indicate that they are
not considered in the objective function and hence not reliable.

\begin{table}
\begin{centering}
\begin{tabular}{|c|c|c|c|c|c|c|}
\hline 
$\left|T\right|$ & ProbSize & IPTime & SolTime & $B^{LB}$ & DisrCost & ResUse\tabularnewline
\hline 
\hline 
35 & 2.0k x 4.4k & 1.1 & 1.3 & 67 & 11.525 & (2+22)\tabularnewline
\hline 
40 & 2.3k x 5.0k & 24 & 27 & 66 & 11.450 & (2+26)\tabularnewline
\hline 
45 & 2.6k x 5.6k & 113 & 163 & 66 & 11.450 & (2+26)\tabularnewline
\hline 
50 & 2.9k x 6.3k & 23 & 742 & 66 & 11.450 & (2+26)\tabularnewline
\hline 
\end{tabular}
\par\end{centering}
\caption{\label{tab:Traffic-disruption-lb}Traffic disruption lower bound for
different planning freedoms}
\end{table}

First we can establish that the lowest possible $B^{LB}=66$ for this
instance. Secondly we see a linear growth in problem size, but a much
larger growth in total solution time (although the integer solution
is found fairly quickly). This shows the difficulty in proving optimality
when having large scheduling freedom and many symmetric solutions.
Finally we note that a low blocking value comes at the price of a
high resource usage.

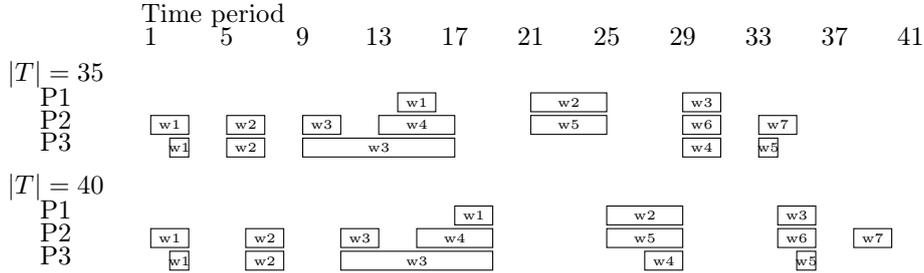
\begin{figure}
\input{figures/ganttsmall_tikz.tex}

\caption{\label{fig:Gantt-block-sol}Schedule solutions with minimal traffic
disruption.}
\end{figure}

The obtained solutions for $\left|T\right|=35$ and $\left|T\right|=40$
are shown as a Gantt view in Figure \ref{fig:Gantt-block-sol}. As
expected, we see that the project tasks are coordinated so as to reduce
the traffic disruption cost. The two solutions have the same structure
-- the only principal difference is that the two last tasks of P3
are placed differently. The latter solution has a lower blocking value,
but cannot be achieved for $\left|T\right|=35$ (due to the minimum
rest time requirement of 8 time periods after P3/W3 in combination
with the task lengths and min rest times for the final tasks of the
projects). This example shows the interplay of schedule length, task
settings, blocking values and cost factors. Finally, we note that
the initial task groups are spaced differently in the two solutions,
which does not affect the objective value. Hence the task spacing
is arbitrary, unless restricted by other factors. Time-dependent blocking
costs, e.g. letting time periods with little traffic have lower blocking
cost, would remedy this arbitrariness (and reduce the symmetries of
the problem).

We now investigate the effect of supplying the valid lower bound $B^{LB}\geq66$
to the MIP solver. Table \ref{tab:Impr-with-valid-blocking-lb} shows
the improvement in solving performance. We see that the time for proving
optimality is reduced substantially, while there is less improvement
in the time for finding integer solutions.

\begin{table}
\begin{centering}
\begin{tabular}{|c|c|c|c|}
\hline 
$\left|T\right|$ & Without lb & With lb & Relative change\tabularnewline
\hline 
\hline 
35 & 1.1 : 1.3 & 0.7 : 0.8 & - 36 : 38 \%\tabularnewline
\hline 
40 & 24 : 27 & 26 : 26 & +8 : -4 \%\tabularnewline
\hline 
45 & 113 : 163 & 79 : 79 & - 30 : 51 \%\tabularnewline
\hline 
50 & 23 : 742 & 28 : 28 & +22 : -96 \%\tabularnewline
\hline 
\end{tabular}
\par\end{centering}
\caption{\label{tab:Impr-with-valid-blocking-lb}Solving performance (IPTime
: SolTime) without and with lower bound on blocking value $B^{LB}$}
\end{table}

\subsection{Resource usage\label{subsec:Resource-usage}}

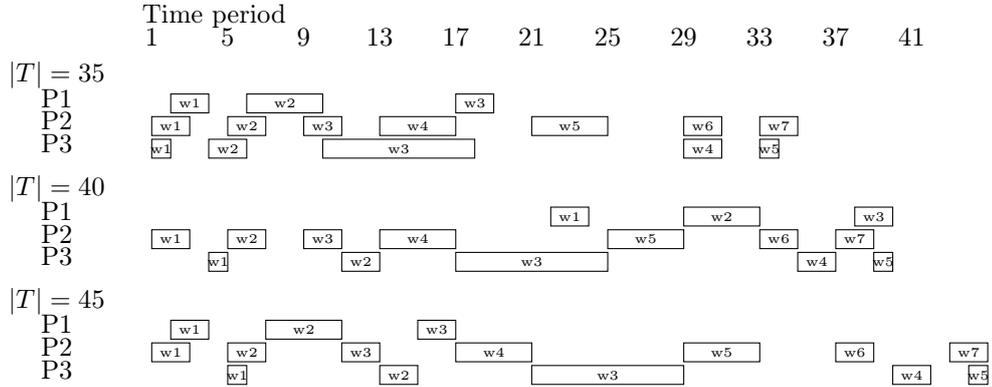
\begin{figure}
\input{figures/gantt_resuse_tikz.tex}

\caption{\label{fig:Gantt-resuse-sol}Schedule solutions with minimal resource
usage.}
\end{figure}

In the next experiment we disregard the traffic disruptions by setting
$c^{disr}=0$, and instead only focus on the resource usage (with
unit resource costs $c_{k}^{K}=1,\forall k$). Again we solve for
the four planning period sizes and obtain the results as shown in
Table \ref{tab:Res-usage-lb} (with the same columns as in Table \ref{tab:Traffic-disruption-lb}).
A Gantt view of the first three solutions are shown in Fig \ref{fig:Gantt-resuse-sol}.
We see that the tasks become separated from each other, which reduces
the resource usage (while increasing the traffic blocking values).
A slight increase in solution time can be noted with increasing planning
period size, but there is no tail effect for proving optimality (SolTime
is very close to IPTime). Thus the solver has no difficulty in establishing
a lower bound for the resource usage (despite the existence of solution
symmetries) as opposed to what was observed for the traffic disruption
cost in Table \ref{tab:Traffic-disruption-lb}.

\begin{table}
\begin{centering}
\begin{tabular}{|c|c|c|c|c|c|c|}
\hline 
$\left|T\right|$ & ProbSize & IPTime & SolTime & $B^{LB}$ & DisrCost & ResUse\tabularnewline
\hline 
\hline 
35 & 2.0k x 4.4k & 0.4 & 0.4 & (130) & (16.250) & 2+18\tabularnewline
\hline 
40 & 2.3k x 5.0k & 18.6 & 18.8 & (180) & (20.000) & 1+12\tabularnewline
\hline 
45 & 2.6k x 5.6k & 22.1 & 22.1 & (185) & (20.375) & 1+10\tabularnewline
\hline 
50 & 2.9k x 6.3k & 34.6 & 34.6 & (185) & (20.375) & 1+10\tabularnewline
\hline 
\end{tabular}
\par\end{centering}
\caption{\label{tab:Res-usage-lb}Resource usage lower bound for different
planning freedoms}
\end{table}

We now check the effect of introducing a simple valid lower bound,
$q_{k}\geq\left(1-s_{u}\right)\max_{w}(g_{uwk})\:\forall u,k$, on
the resource usage. The improvement in solution times is listed in
Table \ref{tab:Impr-with-valid-resuse-lb}. Once again we see an improvement
in solution times---this time both for finding integer solutions
and for proving optimality. We also note that there is no increase
in solution time with increasing planning period size (for this problem
instance).

\begin{table}
\begin{centering}
\begin{tabular}{|c|c|c|c|}
\hline 
$\left|T\right|$ & Without lb & With lb & Relative change\tabularnewline
\hline 
\hline 
35 & 0.4 : 0.4 & 0.6 : 0.6 & +25 : 25 \%\tabularnewline
\hline 
40 & 18.6 : 18.8 & 9.8 : 9.8 & - 47 : 48 \%\tabularnewline
\hline 
45 & 22.1 : 22.1 & 10.7 : 10.7 & - 52 : 52 \%\tabularnewline
\hline 
50 & 34.6 : 34.6 & 8.5 : 8.5 & - 75 : 75 \%\tabularnewline
\hline 
\end{tabular}
\par\end{centering}
\caption{\label{tab:Impr-with-valid-resuse-lb}Solving performance (IPTime
: SolTime) without and with lower bound on resource usage}
\end{table}

\subsection{Combined traffic disruption and resource usage cost\label{subsec:Combined-traffic-and-resources}}

\begin{figure}
\input{figures/gantt_both_tikz.tex}

\caption{\label{fig:Gantt-both-sol}Schedule solutions with optimal traffic
disruption plus resource usage cost}
\end{figure}
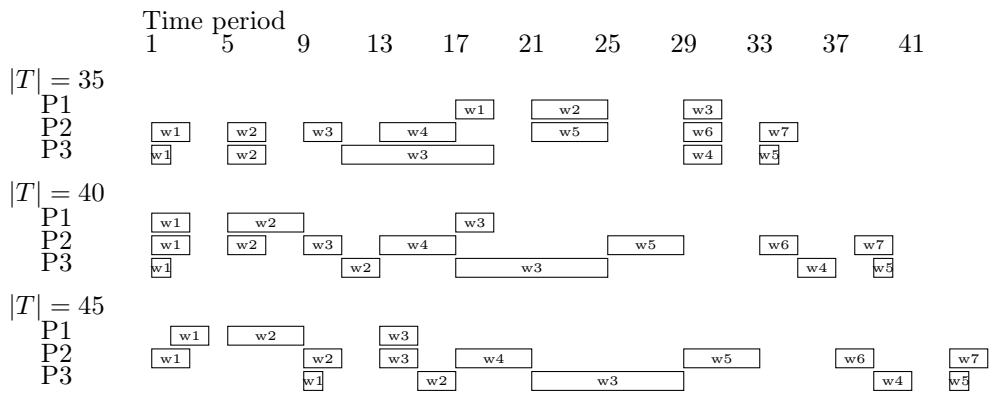

The final experiment concerns the problem of optimizing the sum of
traffic disruption and resource usage cost. We perform this test with
all combinations of supplying valid lower bounds for blocking value
and resource usage. The obtained solutions are presented in Table
\ref{tab:Combined-result} and the performance comparison for the
valid lower bounds are presented in Table \ref{tab:Combined-gen-lb}
and \ref{tab:Combined-spec-lb}, where the first table shows the result
of supplying the general lower bounds ($B^{LB}\geq66$, $q_{R1}\geq1$,
and $q_{R2}\geq10$) while the second table gives the performance
when supplying specific lower bounds for the particular planning problem
size (as obtained from Tables \ref{tab:Traffic-disruption-lb} and
\ref{tab:Res-usage-lb}). Columns two to five list the solving performance
(IPTime : SolTime) in seconds without supplied lower bounds (No lb),
with valid bounds on blocking value (B lb), resource usage (Res lb)
and both blocking value and resource usage (B+Res lb). For reference
we also show the Gantt view of the first three solutions in Figure
\ref{fig:Gantt-both-sol}.

The solutions use the minimal amount of resources while the traffic
disruption cost goes up, which is in accordance with the chosen cost
factors. Hence the work tasks are primarily separated from each other
(so as to reduce resource usage), and only grouped together (for reduced
traffic disruption) as long as it does not increase the resource usage.
The primary performance gain comes from the valid lower bound on resource
usage, while the lower bound on the blocking value has less effect.
Supplying both lower bounds reduces the time for proving optimality
but sometimes hampers the time for finding the best integer solution.
Specific lower bound values further improve the solving performance,
but might not be justified when considering the additional effort
of establishing these specific values.

\begin{table}
\begin{centering}
\begin{tabular}{|c|c|c|c|}
\hline 
$\left|T\right|$ & $B^{LB}$ & DisrCost & ResUse\tabularnewline
\hline 
\hline 
35 & 69 & 11.675 & 2+18\tabularnewline
\hline 
40 & 95 & 13.625 & 1+12\tabularnewline
\hline 
45 & 101 & 14.075 & 1+10\tabularnewline
\hline 
50 & 101 & 14.075 & 1+10\tabularnewline
\hline 
\end{tabular}
\par\end{centering}
\caption{\label{tab:Combined-result}Results for combined traffic disruption
and resource usage cost}
\end{table}

\begin{table}
\begin{centering}
\begin{tabular}{|c|c|c|c|c|}
\hline 
$\left|T\right|$ & No lb & B lb & Res lb & B+Res lb\tabularnewline
\hline 
\hline 
35 & 5.3 : 5.3 & 3.3 : 3.4 & 4.6 : 4.6 & 4.4 : 4.5\tabularnewline
\hline 
40 & 30 : 38 & 39 : 42 & 41 : 45 & 27 : 29\tabularnewline
\hline 
45 & 333 : 393 & 197 : 240 & 124 : 221 & 61 : 206\tabularnewline
\hline 
50 & 161 : 1764 & 110 : 1773 & 304 : 1305 & 418 : 1300\tabularnewline
\hline 
\end{tabular}
\par\end{centering}
\caption{\label{tab:Combined-gen-lb}Performance without and with general lower
bound values}
\end{table}

\begin{table}
\begin{centering}
\begin{tabular}{|c|c|c|c|c|}
\hline 
$\left|T\right|$ & No lb & B lb & Res lb & B+Res lb\tabularnewline
\hline 
\hline 
35 & 5.3 : 5.3 & 4.1 : 4.1 & 3.8 : 3.9 & 0.6 : 2.8\tabularnewline
\hline 
40 & 30 : 38 & (39 : 42) & 32 : 32 & 12 : 14\tabularnewline
\hline 
45 & 333 : 393 & (197 : 240) & 160 : 212 & 33 : 149\tabularnewline
\hline 
50 & 161 : 1764 & (110 : 1773) & 293 : 1302 & 325 : 1130\tabularnewline
\hline 
\end{tabular}
\par\end{centering}
\caption{\label{tab:Combined-spec-lb}Performance without and with specific
lower bound values}
\end{table}

From these computational experiments we first conclude that symmetry
reduction, shorter planning horizons, and smaller scheduling freedom
for the work tasks are the primary factors for achieving shorter solution
times. Supplying valid lower bounds can also improve the solving performance
substantially. In particular it helps to give a lower bound for the
aspect that dominates the objective value---in this case the resource
usage.

\section{Conclusions and discussion\label{sec:Conclusions}}

The bi-level approach which have been described in this report fits
the requirements of renewal project scheduling and traffic flow routing.
In particular, it decomposes the problem into its two major aspects
and enables the use of different mathematical formulations and time
resolutions in the two problem levels. Thus, two important combinatorial
problems can be tackled, one for solving the coordination of work
tasks for several upgrade projects, and the other for optimally adjusting
the network traffic due to imposed capacity restrictions. By iterating
between these two problems, the combined planning problem as described
in Section \ref{subsec:Background-problem-description} can be addressed.

So far, we have only been able to conduct computational experiments
on the two models separately and the integration into one bi-level
model remains to be done. However, the network flow model has so far
shown good performance and at this moment we are positive about the
tractability of this formulation. Furthermore, this flow approach
for handling network wide railway traffic on a macroscopic scale is
useful in its own right, particularly when analyzing transportation
plans, traffic flows and network capacity without the need for any
timetable data.

The mixed integer model for project scheduling has so far only been
evaluated on small problem instances. There are some concerns regarding
the scalability of this approach, and it might be necessary to use
some heuristic or matheuristic approach when treating larger planning
problems. However, the handling of traffic blockings and the lower/upper
bound estimation of traffic impact from concurrent capacity restrictions
enables the upper level to consider the potential traffic disruptions
before invoking the detailed lower level computations. It might also
be possible to calibrate, perhaps in an adaptable approach, the level
of possible traffic coordination (parameter $\zeta$) for different
network, traffic and project relations.

Finally, we have so far not treated the multi-annual planning problem.
In that setting, some additional features need to be introduced---primarily
budget limitations between years (when allowing projects to be moved
between years), and possibly (discounted) values of increasing quality
costs when renewals are delayed vs service improvements after project
completion. In any case, the bi-level approach described in this document,
is a promising starting point for the multi-annual planning case.

\printbibliography[heading=bibintoc]

\pagebreak

\appendix

\section*{Appendices}

\section{Sketch of set partitioning / column generation model\label{app:Set-partition-model}}

The scheduling of upgrade projects and traffic flow could be formulated
as a set partitioning / column generation model. Here we roughly sketch
this, by first considering an integrated approach with two types of
columns: one for the traffic assignment and one for upgrade project
scheduling. Given the major difference in time resolution for these
two aspects we then sketch a bi-level approach, where the upper level
schedules the upgrade projects, while the lower level treats the subsequent
adjustment of the traffic.

A common notation is used as follows. For each traffic relation there
is a limited number of routing options, but an enormous amount of
timed flow assignments that schedule the trains along these routing
options:

\noindent %
\begin{tabular}{c>{\raggedright}p{0.9\textwidth}}
\noalign{\vskip3pt}
$\mathbb{P}_{r}$ & all (considered) timed flow assignments $p$ for traffic relation
$r$. \tabularnewline
\noalign{\vskip3pt}
$\mathbf{a}_{rp}$ & a column vector for flow assignment $p\in\mathbb{P}_{r}$ (of traffic
relation $r$), where the column values $a_{rp}^{lt}$ are the number
of trains that utilize link $l$ in time period $t$.\tabularnewline
\noalign{\vskip3pt}
$x_{rp}$ & binary variable corresponding to each $\mathbf{a}_{rp}$.\tabularnewline
\noalign{\vskip3pt}
$c_{rp}^{T}$ & the attributed cost value for flow assignment $p\in\mathbb{P}_{r}$.\tabularnewline[3pt]
\noalign{\vskip3pt}
\end{tabular}

\noindent For each upgrade project $u\in U$ there are also an enormous
amount of possible work schedules:

\noindent %
\begin{tabular}{c>{\raggedright}p{0.9\textwidth}}
\noalign{\vskip3pt}
$\mathbb{V}_{u}$ & all (considered) work schedules for upgrade project $u$.\tabularnewline
\noalign{\vskip3pt}
$\mathbf{b}_{uv}$ & a column vector for work schedule $v\in\mathbb{V}_{u}$ (of upgrade
$u$), where the column values $b_{uv}^{lt}$ are the number of capacity
blockings caused on link $l$ in time period $t$.\tabularnewline
\noalign{\vskip3pt}
$y_{uv}$ & binary variable corresponding to $\mathbf{b}_{uv}$.\tabularnewline
\noalign{\vskip3pt}
$c_{uv}^{W}$ & the attributed cost value for work schedule $v\in\mathbb{V}_{u}$
for upgrade $u$.\tabularnewline[3pt]
\noalign{\vskip3pt}
\end{tabular}

\subsection{Integrated approach for joint project and traffic scheduling}

In an integrated approach the master problem contains both column
types, and two sub-problems will be used for generating new cost reducing
columns. The master problem is to find flow assignments for all traffic
relations and work schedules for all upgrade projects that minimizes
the objective

\[
\sum_{r\in R,p\in\mathbb{P}_{r}}c_{rp}^{T}x_{rp}+\sum_{u\in U,v\in\mathbb{V}_{u}}c_{uv}^{W}y_{uv}
\]
 subject to the restrictions that:

Exactly one routing configuration $p$ should be allocated to each
service $r$

\[
\forall r\in R:\,\sum_{p\in\mathbb{P}_{r}}x_{rp}=1;
\]

Exactly one schedule $v$ for each upgrade $u$ should be chosen

\[
\forall u\in U:\,\sum_{v\in\mathbb{V}_{u}}y_{uv}=1;
\]

And for each link $l$ and time period $t$, the sum of allocated
trains and the number of capacity blockings must not exceed the nominal
link capacity

\[
\forall l\in L,t\in T:\,\sum_{r\in R,p\in\mathbb{\mathbb{P}}_{r}}a_{rp}^{lt}x_{rp}+\sum_{u\in U,v\in\mathbb{V}_{u}}b_{uv}^{lt}y_{uv}\leq C_{l}^{nom}.
\]

\noindent The sub-problem pricing procedures will then use the dual
values from these constraints to generate potential cost reducing
columns. However, as noted above, the time horizon as well as the
time resolution of these columns differ substantially.

\subsection{Bi-level approach}

The basis for the bi-level approach is to split the two parts of the
problem into two different solving processes. The upper level is the
scheduling of the upgrade projects according to an initial cost function.
For each operational day which is affected by one or more (temporary)
capacity restriction (TCR), the subsequent traffic adjustments are
calculated by the lower level problem. Thus, only the affected traffic
days needs to be considered in the lower level, which reduces the
problem size. The obtained traffic adjustments can then be considered
by the upper level problem in order to find a better TCR schedule,
taking both the cost for the project schedule and the traffic flow
into account. The two processes are iterated until some convergence
or optimality criteria is met.

This approach opens up for using different objective functions in
the two levels, as well as using different modeling and solution approaches.
Furthermore, we can use different discretisations of time on the upper
and lower level, which suits the two different types of schedules.
For upgrades, the schedules are on macro level, i.e. weeks and days,
possibly down to hours. For transports, the schedules are on a finer
detail, from hours down to minutes.

We now introduce some further notation, so as to be able to consider
resource usage in the upper level, and to only consider disturbed
days in the lower level, as follows:

\noindent %
\begin{tabular}{cl}
\noalign{\vskip3pt}
$WD(\mathbf{y})$ & set of disrupted (working) days for upgrade project schedule $\mathbf{y}$\tabularnewline
\noalign{\vskip3pt}
$F_{d}(\mathbf{x}^{d})$ & disruption cost for traffic assignment $\mathbf{x}^{d}$ of day $d$\tabularnewline
\noalign{\vskip3pt}
$K$ & the set of all resources\tabularnewline
\noalign{\vskip3pt}
$Q_{k}$ & number of available resources $k$\tabularnewline
\noalign{\vskip3pt}
$L_{k}^{K}$ & set of links that resource $k$ can service\tabularnewline
\noalign{\vskip3pt}
$T_{k}^{K}$ & set of time periods that resource $k$ can service\tabularnewline
\noalign{\vskip3pt}
$q_{k}$ & (variable) max resource usage of resource $k$ in any work period
$t$\tabularnewline
\noalign{\vskip3pt}
$g_{uv}^{ktl}$ & usage of resource $k$ during time instance $t$ on link $l$, for
upgrade $u$ and schedule $v$\tabularnewline[3pt]
\noalign{\vskip3pt}
\end{tabular}

\subsubsection{The upper level}

The upper level problem schedules the upgrade projects. The objective
is to minimize the sum of work costs, resource usage, and traffic
penalties based on the adjusted traffic from the lower level.

\[
\text{Min}\sum_{u\in U,v\in V}c_{uv}^{W}y_{uv}+\sum_{k\in K}c_{k}^{K}q_{k}+\sum_{d\in WD(\mathbf{y})}F_{d}(\mathbf{x}^{d})
\]
subject to restrictions that:

Exactly one schedule $v$ should be chosen for each upgrade $u$

\[
\forall u\in U:\sum_{v\in\mathbb{V}_{u}}y_{uv}=1
\]

Consumed capacity on link $l$ for time period $t$ must be respected
(where $C_{lt}$ is the estimated traffic capacity need)

\[
\forall t\in T,l\in L:\,\sum_{u\in U,v\in V_{u}}b_{uv}^{lt}y_{uv}\leq C_{l}^{nom}-C_{lt}
\]

Usage of resource $k$ in any time period cannot exceed $Q_{k}$

\[
\forall k\in K,t\in T_{k}^{K}:\,\sum_{l\in L_{k}^{K}}g_{uv}^{ktl}\sum_{u\in U,v\in v_{k}}y_{uv}\leq q_{k}\leq Q_{k}
\]
Note that the values $C_{lt}$ constitute the capacity connection
between the upper and lower level problems.

\subsubsection{The lower level}

The lower level problem schedules the trains, with the capacity restrictions
from the upper level as input. The objective is to find an optimal
traffic schedule. Note that this is not necessarily the schedule that
through-puts as much trains as possible, other factors may also affect
what is an optimal schedule such as conformity over the year, regularity,
robustness as well as steadiness.

We will first concentrate on the individual trains and a cyclic period,
where the objective is to find assignments of transport tasks, defined
by origin-destination pairs with periodic repetition.

\[
\text{Min}\sum_{r\in R,p\in Pr}c_{rp}^{T}x_{rp}
\]
subject to restrictions that:

Exactly one flow assignment (routing and timing) $p$ is chosen for
each transport $r$

\[
\forall r\in R:\sum_{p\in P_{r}}x_{rp}=1
\]

Consumed capacity by traffic must respect the capacity restrictions
$C_{lt}$ on each link $l$ and time period $t$, which is the remaining
capacity according to the current upgrade project schedule $\mathbf{y}$

\[
\forall t\in T,l\in L:\,\sum_{r\in R,p\in P_{r}}a_{rp}^{lt}x_{rp}\leq C_{lt}
\]

By letting the traffic schedules be either calculated as cyclic with
repetition or by solving instances of consecutive time frames, longer
periods can be calculated in the lower level model. In both cases,
it will be possible to calculate the cost $F_{d}(\mathbf{x}^{d})$
for each working day $d$.

\section{Flow based traffic assignment idea\label{app:Flow-models}}

The following is a write-up of two flow based formulations which were
done in parallel with implementation and testing. The first attempt
is a model that only distinguishes between train types. After that
we describe a multi-commodity model, which correctly will handle the
reachability requirements of all origin-destination relations, but
allows for flow splitting and merging in the nodes. The two final
subsections discuss the problem sizes of these two formulations and
reports on the initial computational experiments.

\subsection{Aggregating over train types}

The basic notation for planning period, railway network and wanted
traffic is the same as described in Section \ref{subsec:Basic-notation}.
The traffic assignment problem can be represented by a directed graph
$\mathcal{G}\left(\mathcal{V},\mathcal{A}\right)$ with a vertex set
$\mathcal{V}\coloneqq\left\{ v_{nth}|n,t,h\in N\times T\times H\right\} $
and an arc set $\mathcal{A}\coloneqq\mathcal{A}^{Tr}\cup\mathcal{A}^{Ni}$
which is the union of link transportation arcs $\mathcal{A}^{Tr}$
and node inventory arcs $\mathcal{A}^{Ni}$.

Each vertex $v_{nth}\in\mathcal{V}$ collects the number of trains
of type $h$ that start, pass or end in node $n$ during time period
$t$. Let $n(i)$, $t(i)$, and $h(i)$ mean the node, the time period
and train type of vertex $i$ respectively. Then the link transportation
arcs $a_{ij}\in\mathcal{A}^{Tr}$ describe all possible train traffic
from vertex $i$ to vertex $j$ such that $\left(n(i),n(j)\right)\in L$,
$t(i)\leq t(j)\leq t(i)+1$, and $h(i)=h(j)$. Similarly, the node
inventory arcs $a_{ij}\in\mathcal{A}^{Ni}$ describe possible train
dwelling at stations between vertex $i$ and $j$ such that $n(i)=n(j)$,
$t(i)+1=t(j)$, and $h(i)=h(j)$. Thus the graph will be partitioned
according to the train types and transportation can only move forward
in time, but must also arrive no later than in the next time period.
Furthermore we impose the restriction that the travel time for any
train type along link $l$ must be less than or equal to $\min\left(\delta_{t}\right)/2$,
i.e. the shortest time period length divided by two. The consequence
of this restriction is that at least half the number of trains can
arrive in the same time period under the assumption that the trains
depart evenly distributed over the time period. The net effect is
that time period lengths must be chosen sufficiently large so as to
fulfill this restriction.

We now define the following set of variables:

\begin{tabular}{l>{\raggedright}p{0.6\columnwidth}}
\noalign{\vskip3pt}
$x_{r}^{C}\in\mathbb{R}_{+}^{R}$ & Number of canceled trains for traffic relation $r$\tabularnewline
\noalign{\vskip3pt}
$x_{rt}^{O},x_{rt}^{D}\in\mathbb{R}_{+}^{R\times T}$ & Number of trains departing (from the origin), arriving (to the destination)
in time period $t$ for traffic relation $r$\tabularnewline
\noalign{\vskip3pt}
$x_{i}^{OD}\in\mathbb{R}^{N\times T\times H}$ & The aggregated source/sink balance of vertex $i$ (auxilliary convenience
variable)\tabularnewline
\noalign{\vskip3pt}
$x_{ij}\in\mathbb{R}_{+}^{(2L+N)\times T\times H}$ & Number of trains along arc $a_{ij}\in\mathcal{A}$\tabularnewline[3pt]
\noalign{\vskip3pt}
\end{tabular}

\noindent and introduce the following convenience notation

\begin{tabular}{l>{\raggedright}p{0.7\columnwidth}}
\noalign{\vskip3pt}
$T_{r}^{O},T_{r}^{D}\subseteq T$ & The allowed departure / arrival time periods (at the origin / destination
nodes respectively) for traffic relation $r$\tabularnewline
\noalign{\vskip3pt}
$\mathcal{A}_{l}^{Tr}\subset\mathcal{A}^{Tr}$ & All link transportation arcs along link $l$\tabularnewline
\noalign{\vskip3pt}
$\mathcal{A}_{lh}^{Tr}\subset\mathcal{A}_{l}^{Tr}$ & All link transportation arcs along link $l$ for train type $h$\tabularnewline[3pt]
\noalign{\vskip3pt}
\end{tabular}

Using this notation we now formulate the fundamental constraints that
a traffic assignment should fulfill. First, the source and sink flows
of running trains should balance the demand minus the number of canceled
trains
\begin{equation}
\forall r\in R:\sum_{t\in T_{r}^{O}}x_{rt}^{O}=\sum_{t\in T_{r}^{D}}x_{rt}^{D}=\sigma_{r}-x_{r}^{C}.
\end{equation}
The source/sink flows of each vertex is then aggregated from all relations
\[
\forall v_{nth}\in\mathcal{V}:x_{i}^{OD}=\sum_{r\in R:o_{r}=n,h_{r}=h}x_{rt}^{O}-\sum_{r\in R:d_{r}=n,h_{r}=h}x_{rt}^{D}
\]
and used in the flow balance constraint
\begin{equation}
\forall v_{nth}\in\mathcal{V}:x_{j}^{OD}+\sum_{a_{ij}\in\mathcal{A}}x_{ij}-\sum_{a_{jk}\in\mathcal{A}}x_{jk}=0\label{eq:flow-balance}
\end{equation}
A limitation is imposed on how many trains that can traverse a link
within the same time period $t$ (with a slight abuse of notation)

\begin{multline}
\forall v_{nth},v_{mth}\in\mathcal{V}:\left(n,m\right)=l\in L:\\
x_{nm}\leq\sum_{a_{in}\in\mathcal{A}:t(i)<t}x_{in}+\left(1-\nu_{lth}\right)x_{n}^{OD}\\
+\sum_{a_{in}\in\mathcal{A}_{ph}^{Tr}:t(i)=t}\frac{1-\nu_{pth}-\nu_{lth}}{1-\nu_{pth}}x_{in}\label{eq:share-within-period}
\end{multline}
where the parameter value $\nu_{lth}$ is the share of time period
$t$ that a train of type $h$ will take to traverse link $l$, which
we use as an estimate for the share of trains departing in time period
$t$ that will arrive in time period $t+1$. The three terms in the
right hand side expression capture the available flow until the point
in time when departure must be made so as to arrive within the same
time period. The terms are: (1) incoming trains (both transportation
and inventory) from the previous time period, (2) available source/sink
balance (under the assumption that these flows are evenly distributed
over the time period), and (3) incoming transportation flow from preceding
links $p$ that departed in the same time period. In the last term,
the fraction expresses the amount of flow from the previous link $p$
that can traverse the successor $l$ within the same time period.\footnote{The computational experiments revealed that constraint (\ref{eq:share-within-period})
does not work so well for flow splitting (but works for flow merging).
In the multi-commodity model we give some revised constraint formulations
to better account for this.}

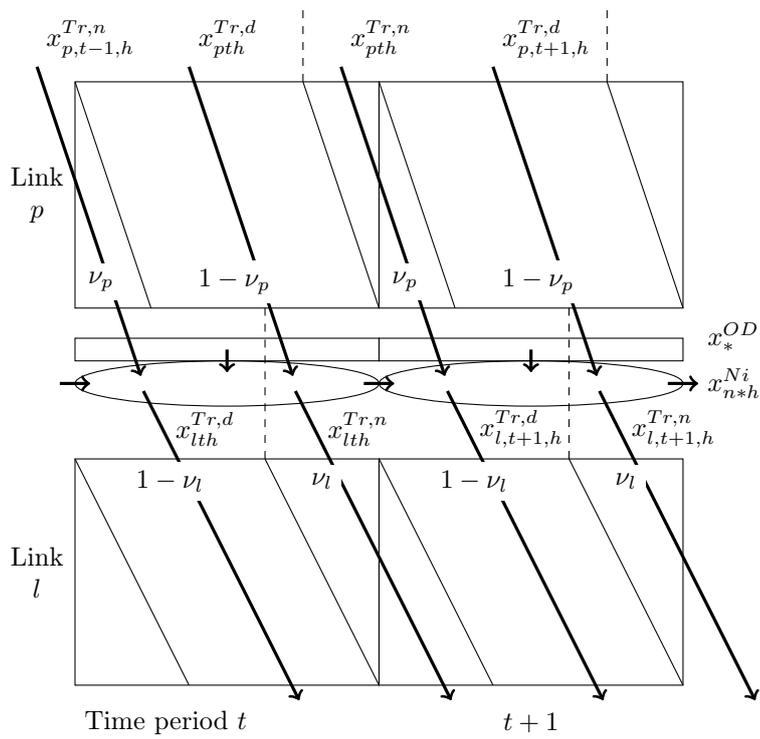
\begin{figure}
\input{figures/flowshare_tikz.tex}

\caption{\label{fig:flow-share}Illustration of flow volume relations. Incoming
transportation flow from link $p$ and outgoing transportation flow
on link $l$.}
\end{figure}

Figure \ref{fig:flow-share} illustrates the different terms and the
flow volume relations. In this figure we use the following convenience
notation: $x_{lth}^{Tr,d}$ (transportation flow along link $l$,
starting and arriving in time period $t$, for train type $h$), $x_{lth}^{Tr,n}$
(transportation flow along link $l$ between time period $t$ and
$t+1$ for train type $h$), and $x_{nth}^{Ni}$ (node inventory flow
for node $n$ and train type $h$ between time period $t$ and $t+1$).
The dashed lines mark the cut-off times for flow that can arrive within
the same time period, which determine the factors used in the second
and third right hand side terms of constraint (\ref{eq:share-within-period}).

The line capacity (total number of trains over a link during a time
period) must be respected
\begin{equation}
\forall l\in L,t\in T:\frac{1}{2}\left[\sum_{a_{ij}\in\mathcal{A}_{l}^{Tr}:t(i)=t}x_{ij}+\sum_{a_{ij}\in\mathcal{A}_{l}^{Tr}:t(j)=t}x_{ij}\right]\leq C_{lt}.\label{eq:capacity}
\end{equation}

\noindent where departing (first term) and arriving (second term)
trains will acquire one half each of the available capacity.

Note that constraint (\ref{eq:share-within-period}) can be made direction
dependent if trains have different traversal times along the up/down
link direction. The model can also be made cyclic by adding wrap around
node inventory arcs -- possibly at the risk of hampering solving
performance. Finally it is perfectly possible to make the capacity
constraint (\ref{eq:capacity}) direction dependent, or to use both
directional and total capacities as in \textcite{liden_optimization_2017}.

Merely having flow balance does however not guarantee reachability
for all relations. In fact, the above constraints will allow for train
flow to arrive before it has departed. Hence we add constraints for
respecting the minimum possible travel time

\[
\forall r\in R,t\in T:\sum_{t'=1}^{t}x_{rt'}^{D}\leq\sum_{t'=1}^{\left\lfloor t+1-mt_{r}\right\rfloor }\frac{\min\left(\delta_{t'},\beta_{t}+\delta_{t}-mt_{r}-\beta_{t'}\right)}{\delta_{t'}}x_{rt'}^{O}
\]

Still it is possible to circumvent the reachability requirements for
the traffic relations by mixing flows from different origins. To some
extent it will help to penalize node inventory arcs (or even removing
them entirely), but the model can still overcome capacity restrictions
by flow mixing. This property can be considered as a feature to resolve
blockages, since it allows for short-turning trains. For relations
requiring reachability (like freight trains) it might be possible
to enforce this property with constraints (possibly dynamic) of the
following structure
\[
\forall r\in R^{reach},L_{r}^{trav}:\sum_{a_{ij}\in\bigcup_{L_{r}^{trav}}A_{lh}^{Tr}}x_{ij}\geq\sum_{r\in R^{reach}}\sum_{t\in T}x_{rt}^{O}
\]
where $R^{reach}$ are all relations with reachability requirements
and $L_{r}^{trav}$ is a subset of links (along a certain direction)
which must be traversed by a relation $r$ with reachability requirements.
Unfortunately, computational experiments indicate that the amount
of flow mixing will be complicated to handle and it has been hard
to formulate constraints that resolve this issue.

\subsection{Multi-commodity model}

We now describe the corresponding multi-commodity approach (over all
relations) and compare the two models regarding problem size.

In this formulation the directed graph $\mathcal{G}\left(\mathcal{V},\mathcal{A}\right)$
has a vertex set $\mathcal{V}\coloneqq\left\{ v_{ntr}|n,t,r\in N\times T\times R\right\} $
where the vertices are partitioned according to the traffic relations
$R$. The arc set $\mathcal{A}\coloneqq\mathcal{A}^{Tr}\cup\mathcal{A}^{Ni}$
is the union of link transportation arcs $\mathcal{A}^{Tr}$ and node
inventory arcs $\mathcal{A}^{Ni}$. Each vertex $v_{ntr}\in\mathcal{V}$
collects the number of trains of traffic relation $r$ that start,
pass or end in node $n$ during time period $t$. Let $n(i)$, $t(i)$,
and $r(i)$ mean the node, the time period and traffic relation of
vertex $i$ respectively. Then the link transportation arcs $a_{ij}\in\mathcal{A}^{Tr}$
describe all possible train traffic from vertex $i$ to vertex $j$
such that $\left(n(i),n(j)\right)\in L$, $t(i)\leq t(j)\leq t(i)+1$,
and $r(i)=r(j)$. Similarly, the node inventory arcs $a_{ij}\in\mathcal{A}^{Ni}$
describe possible train dwelling at stations between vertex $i$ and
$j$ such that $n(i)=n(j)\in N$, $t(i)+1=t(j)$, and $r(i)=r(j)$.
We use the same assumptions as before, namely that trains must move
forward in time, that they will arrive no later than in the next time
period, and that at most half of them can arrive within the same time
period.

The following set of variables are used:

\begin{tabular}{l>{\raggedright}p{0.6\columnwidth}}
\noalign{\vskip3pt}
$x_{r}^{C}\in\mathbb{R}_{+}^{R}$ & Number of canceled trains for traffic relation $r$\tabularnewline
\noalign{\vskip3pt}
$x_{rt}^{O},x_{rt}^{D}\in\mathbb{R}_{+}^{R\times T}$ & Number of trains departing (from the origin), arriving (to the destination)
in time period $t$ for traffic relation $r$\tabularnewline
\noalign{\vskip3pt}
$x_{i}^{OD}\in\mathbb{R}^{N\times T\times R}$ & The aggregated source/sink balance of vertex $i$ (auxilliary convenience
variable)\tabularnewline
\noalign{\vskip3pt}
$x_{ij}\in\mathbb{R}_{+}^{(2L+N)\times T\times R}$ & Number of trains along arc $a_{ij}\in\mathcal{A}$\tabularnewline[3pt]
\noalign{\vskip3pt}
\end{tabular}

\noindent along with the convenience notation

\begin{tabular}{l>{\raggedright}p{0.7\columnwidth}}
\noalign{\vskip3pt}
$T_{r}^{O},T_{r}^{D}\subseteq T$ & The allowed departure / arrival time periods (at the origin / destination
nodes respectively) for traffic relation $r$\tabularnewline
\noalign{\vskip3pt}
$\mathcal{A}_{l}^{Tr}\subset\mathcal{A}^{Tr}$ & All link transportation arcs along link $l$\tabularnewline
\noalign{\vskip3pt}
$\mathcal{A}_{lr}^{Tr}\subset\mathcal{A}_{l}^{Tr}$ & All link transportation arcs along link $l$ for traffic relation
$r$\tabularnewline[3pt]
\noalign{\vskip3pt}
\end{tabular}

In order to limit the size of the model we introduce

\begin{tabular}{c>{\raggedright}p{0.75\columnwidth}}
\noalign{\vskip3pt}
$L_{r}^{P}\subseteq L$ & All links that are possible to use for traffic relation $r\in R$\tabularnewline
\noalign{\vskip3pt}
$N_{r}^{P}\subseteq N$ & All nodes that are possible to visit for traffic relation $r\in R$,
as derived from $L_{r}^{P}$\tabularnewline[3pt]
\noalign{\vskip3pt}
\end{tabular}

Using this notation we now reformulate the necessary constraints for
a valid traffic assignment. The source and sink flows of running trains
are the same as before
\begin{equation}
\forall r\in R:\sum_{t\in T_{r}^{O}}x_{rt}^{O}=\sum_{t\in T_{r}^{D}}x_{rt}^{D}=\sigma_{r}-x_{r}^{C}.
\end{equation}

\noindent The source/sink flows of each vertex is collected
\[
\forall r\in R,n\in N_{r}^{P},v_{ntr}\in\mathcal{V}:x_{i}^{OD}=\begin{cases}
x_{rt}^{O} & n=o_{r}\\
0 & n\notin\left\{ o_{r},d_{r}\right\} \\
-x_{rt}^{D} & n=d_{r}
\end{cases}
\]
 and used in the flow balance constraint
\begin{equation}
\forall r\in R,n\in N_{r}^{P},v_{ntr}\in\mathcal{V}:x_{j}^{OD}+\sum_{a_{ij}\in\mathcal{A}}x_{ij}-\sum_{a_{jk}\in\mathcal{A}}x_{jk}=0\label{eq:flow-balance-mc}
\end{equation}

\noindent The limitation on how many trains that can traverse a link
within the same time period $t$

\begin{multline}
\forall r\in R;l\in L_{r}^{P};v_{ntr},v_{mtr}\in\mathcal{V}:\left(n,m\right)=l:\\
x_{nm}\leq\sum_{a_{in}\in\mathcal{A}:t(i)<t}x_{in}+\left(1-\nu_{lth}\right)x_{n}^{OD}\\
+\sum_{a_{in}\in\mathcal{A}_{pr}^{Tr}:t(i)=t}\frac{1-\nu_{pth}-\nu_{lth}}{1-\nu_{pth}}x_{in}\label{eq:share-within-period-mc}
\end{multline}
 where we still use parameter values $\nu_{lth}$ for train of type
$h$ (of relation $r$) as before. As noted in the previous section
the above constraint does not work so well for flow splitting. The
following alternative formulation addresses this issue:

\begin{multline*}
\forall r\in R;p\in L_{r}^{P};v_{itr},v_{ntr}\in\mathcal{V}:\left(i,n\right)=p:\\
\sum_{a_{nm}\in\mathcal{A}_{lr}^{Tr}:t(m)=t}x_{nm}\leq\sum_{a_{jn}\in\mathcal{A}:t(j)<t}x_{jn}+\left(1-\nu_{lth}^{*}\right)x_{n}^{OD}\\
+\frac{1-\nu_{pth}-\nu_{lth}^{*}}{1-\nu_{pth}}x_{in}
\end{multline*}
 where $\nu_{lth}^{*}$ is the largest value among the outgoing arcs
of the left hand side summation ($A_{lr}^{Tr}:t(m)=t$).

Based on the above we can formulate a combined version that will handle
all merge/split cases. First we define the largest link traversal
times for all nodes with outgoing arcs, as follows:

\begin{align*}
\forall r\in R,n\in N_{r}^{P}\setminus d_{r}: & \nu_{nth}^{*}=\max\left[\nu_{lth}:(n,m)=l\in L_{r}^{P}\right]
\end{align*}

\noindent Using these values, we can express the limitation on the
flow within the same time period as follows

\begin{multline*}
\forall r\in R,n\in N_{r}^{P}\setminus d_{r},v_{ntr}\in\mathcal{V}:\\
\sum_{a_{nj}\in\mathcal{A}_{lr}^{Tr}:t(j)=t}\frac{1-\nu_{nth}^{*}}{1-\nu_{lth}}x_{ij}\leq\sum_{a_{in}\in\mathcal{A}:t(i)<t}x_{in}+\left(1-\nu_{nth}^{*}\right)x_{n}^{OD}\\
+\sum_{a_{in}\in\mathcal{A}_{pr}^{Tr}:t(i)=t}\frac{1-\nu_{pth}-\nu_{nth}^{*}}{1-\nu_{pth}}x_{in}
\end{multline*}

Finally, the line capacity constraint is the same as previously
\begin{equation}
\forall l\in L,t\in T:\frac{1}{2}\left[\sum_{a_{ij}\in\mathcal{A}_{l}^{Tr}:t(i)=t}x_{ij}+\sum_{a_{ij}\in\mathcal{A}_{l}^{Tr}:t(j)=t}x_{ij}\right]\leq C_{lt}.
\end{equation}

This formulation will assert that all departing trains will arrive
at their destination and subsequently we do not need any arrival restrictions
as in the train type formulation.

\subsection{Problem size\label{subsec:Problem-size}}

The dominating number of variables are the arc flow variables $x_{ij}$,
while the number of constraints are dominated by the flow balance
constraints, (\ref{eq:flow-balance}) and (\ref{eq:flow-balance-mc})
respectively, and the flow share within period constraints, (\ref{eq:share-within-period})
and (\ref{eq:share-within-period-mc}) respectively. Table \ref{tab:Size-estimates}
gives an estimate of the number of variables and constraints for the
two model variants, where letters L, T, H, R and S denote the number
of links, time periods, train types, relations, and average number
of possible links per relation (given by the input data sets $L_{r}^{P}$).
We assume that the network is relatively sparse, such that number
of nodes are in the same order of magnitude as the number of links,
e.g. $\left|N\right|\approx\left|L\right|$.

\begin{table}
\begin{centering}
\begin{tabular}{|c|c|c|}
\hline 
 & Aggregated train types & Multi-commodity\tabularnewline
\hline 
\hline 
Variables & $\approx\textrm{3LTH}$ & $\approx\textrm{3RTS}$\tabularnewline
\hline 
Constraints & $\approx\textrm{2LTH}$ & $\approx\textrm{2RTS}$\tabularnewline
\hline 
\end{tabular}
\par\end{centering}
\caption{\label{tab:Size-estimates}Estimate of number of variables and constraints}
\end{table}

The national railway system of Sweden is divided into a bit more than
200 track links, where all entering traffic must proceed to the end
of the link (although there might be several intermediate stations
for passenger exchange and meet/pass handling). Thus at least 400
uni-directional links will be necessary when treating the complete
railway network. If we assume that one day of traffic is divided into
6 time periods (of 4 hours each) and that we want to distinguish between
5 train types, then a daily problem instance for the aggregated model
will need about $3\times400\times6\times5=36k$ variables and $24k$
constraints. A weekly problem instance will require at least $250k$
variables and $170k$ constraints.

The corresponding problem sizes for the multi-commodity formulation
depends mainly on how many links that each traffic relation can traverse
(letter S), which should include all rerouting possibilities. For
the Swedish case we typically have about 1500 traffic relations, and
if we assume that they on average can traverse 10 track links we get
$3\times1500\times6\times10=270k$ variables and $180k$ constraints
for a daily instance. A weekly problem would require about $1900k$
variables and $1300k$ constraints. If on the other hand, each traffic
relation can traverse 50 or 100 track links, then the problem sizes
will be 5 and 10 times higher.

While we see that the multi-commodity formulation gives a substantially
larger problem, it still seems tractable to handle daily instances
for the complete national network of Sweden. Weekly problems might
be solvable for a modern LP solver, if the possible number of links
per traffic relation is sufficiently limited.

\subsection{Computational experiments}

\begin{figure}
\begin{centering}
\includegraphics[viewport=0bp 0bp 346bp 349bp,scale=0.6]{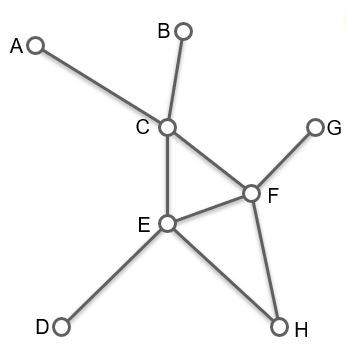}
\par\end{centering}
\caption{\label{fig:Small-flow-case}Small data instance}
\end{figure}

Both variants of the flow-based traffic assignment model have been
implemented in MiniZinc and solved with COIN-BC as MIP solver. A small
data instance have been used, with 8 nodes and 9 uni-directional links,
a scheduling period of 12 time periods and 4 one-way traffic relations.
The network layout is illustrated in Figure \ref{fig:Small-flow-case}
and the properties of the traffic relations are given in Table \ref{tab:Traffic-small-case}.
The last column lists the set of possible links used in the multi-commodity
formulation. The train speed / traversal time share values $\nu_{l*h}$
and the capacity restrictions $C_{l*}$ have been set as listed in
Table \ref{tab:Train-restrictions}, which makes the route C-E-H faster
than C-F-H. A second instance with return traffic included (which
results in 18 links and 8 traffic relations) have also been tested.
Finally, the objective function used for this experiment minimizes
the sum of canceled trains, node inventory flows and the average traveling
time of the traffic relations.

\begin{table}
\begin{centering}
\begin{tabular}{|c|>{\centering}p{1.2cm}|>{\centering}p{1cm}|>{\centering}p{1.5cm}|c|}
\hline 
From/To & Train type & Num. trains & Start/End period & Possible links, $L_{r}^{P}$\tabularnewline
\hline 
\hline 
A/H & freight & 10 & 1/10 & \{A-C, C-E, E-H, C-F, F-H\}\tabularnewline
\hline 
B/H & pax & 25 & 3/11 & \{B-C, C-E, E-H, C-F, F-H\}\tabularnewline
\hline 
D/G & freight & 21 & 1/11 & \{D-E, E-F, F-G\}\tabularnewline
\hline 
E/F & pax & 19 & 3/11 & \{E-F\}\tabularnewline
\hline 
\end{tabular}
\par\end{centering}
\caption{\label{tab:Traffic-small-case}Traffic relations, small data instances}
\end{table}
\begin{table}
\begin{centering}
\begin{tabular}{|c|c|c|c|}
\hline 
Link & $\nu_{l*,fr}$ & $\nu_{l*,px}$ & $C_{l*}$\tabularnewline
\hline 
\hline 
A-C & 0.3 & 0.2 & 5\tabularnewline
\hline 
B-C & 0.2 & 0.1 & 5\tabularnewline
\hline 
C-E & 0.2 & 0.1 & 5\tabularnewline
\hline 
C-F & 0.2 & 0.1 & 5\tabularnewline
\hline 
E-H & 0.3 & 0.2 & 5\tabularnewline
\hline 
F-H & 0.4 & 0.3 & 5\tabularnewline
\hline 
D-E & 0.3 & 0.2 & 5\tabularnewline
\hline 
E-F & 0.2 & 0.1 & 5\tabularnewline
\hline 
F-G & 0.2 & 0.1 & 5\tabularnewline
\hline 
\end{tabular}
\par\end{centering}
\caption{\label{tab:Train-restrictions}Train speed values and capacity settings}
\end{table}
The aggregated train type variant exhibits substantial amount of flow
mixing between the traffic relations (of the same train types) and
we think that overcoming these issues by, for example, introduction
of cuts or lazy constraints is not worth further investigation.

The multi-commodity variant solves these small instances correctly,
but results in a larger problem size (see discussion in Section \ref{subsec:Problem-size}).
For these instances the solving time is negligible (<\textcompwordmark <
1 s) and with no noticeable difference in performance between the
two formulations.

The primary result of these tests are that the method for calculating
the share of flow within the same time period works reasonably well
except for an important issue regarding spreading of fractional flow.
In Table \ref{tab:Flow-solution} the result for one of the traffic
relations is presented, which shows the amount of flow that moves
forward in time. We also note that for this small case the capacity
restrictions are limiting (for links B-C, C-E, and E-H during time
periods 6--8), as shown in Table \ref{tab:Capacity-usage}, and that
some train cancellations are necessary in order to find a feasible
solution.

\begin{table}
\begin{centering}
\begin{tabular}{|c|cccccccccc|}
\hline 
 & \multicolumn{10}{c|}{Time period}\tabularnewline
 & 3 & 4 & 5 & 6 & 7 & 8 & 9 & 10 & 11 & 12\tabularnewline
\hline 
\hline 
$x_{rt}^{O}$ & {\footnotesize{}1} & {\footnotesize{}1} & {\footnotesize{}3} & {\footnotesize{}5} & {\footnotesize{}5} & {\footnotesize{}5} & {\footnotesize{}3} & {\footnotesize{}1} &  & \tabularnewline
\hline 
B-C & {\footnotesize{}0.9:0.1} & {\footnotesize{}0.9:0.1} & {\footnotesize{}2.7:0.3} & {\footnotesize{}4.5:0.5} & {\footnotesize{}4.5:0.5} & {\footnotesize{}4.5:0.5} & {\footnotesize{}2.7:0.3} & {\footnotesize{}0.9:0.1} &  & \tabularnewline
C-E & {\footnotesize{}0.8:0.1} & {\footnotesize{}0.9:0.1} & {\footnotesize{}2.5:--} & {\footnotesize{}4.2:--} & {\footnotesize{}4.2:--} & {\footnotesize{}4.2:--} & {\footnotesize{}2.9:--} & {\footnotesize{}1.1:0.1} & {\footnotesize{}0.1:--} & \tabularnewline
E-H & {\footnotesize{}0.6:0.2} & {\footnotesize{}0.8:0.2} & {\footnotesize{}2.1:0.6} & {\footnotesize{}3.3:0.9} & {\footnotesize{}3.3:0.9} & {\footnotesize{}3.3:0.9} & {\footnotesize{}2.3:0.6} & {\footnotesize{}0.9:0.2} & {\footnotesize{}0.2:0.0} & \tabularnewline
C-F & {\footnotesize{}--:--} & {\footnotesize{}--:--} & {\footnotesize{}--:0.3} & {\footnotesize{}0.1:0.5} & {\footnotesize{}0.3:0.5} & {\footnotesize{}0.3:0.5} & {\footnotesize{}--:0.3} & {\footnotesize{}--:--} &  & \tabularnewline
F-H & {\footnotesize{}--:--} & {\footnotesize{}--:--} & {\footnotesize{}--:--} & {\footnotesize{}0.4:0.0} & {\footnotesize{}0.7:0.1} & {\footnotesize{}0.7:0.1} & {\footnotesize{}0.5:--} & {\footnotesize{}0.3:--} &  & \tabularnewline
\hline 
$x_{rt}^{D}$ & {\footnotesize{}0.6} & {\footnotesize{}1.0} & {\footnotesize{}2.2} & {\footnotesize{}4.2} & {\footnotesize{}4.9} & {\footnotesize{}5} & {\footnotesize{}3.8} & {\footnotesize{}1.8} & {\footnotesize{}0.4} & {\footnotesize{}0.0}\tabularnewline
\hline 
\end{tabular}
\par\end{centering}
\caption{\label{tab:Flow-solution}Flow solution for traffic relation B/H (pax).
Link flow values are given as $x_{ltr}^{Tr,d}$:$x_{ltr}^{Tr,n}$.}
\end{table}

\begin{table}
\begin{centering}
\begin{tabular}{|c|cccccccccc|}
\hline 
 & \multicolumn{10}{c|}{Time period}\tabularnewline
Link & 3 & 4 & 5 & 6 & 7 & 8 & 9 & 10 & 11 & 12\tabularnewline
\hline 
\hline 
B-C & {\footnotesize{}0.95} & {\footnotesize{}1.0} & {\footnotesize{}2.9} & {\footnotesize{}4.9} & \textcolor{red}{\footnotesize{}5.0} & \textcolor{red}{\footnotesize{}5.0} & {\footnotesize{}3.1} & {\footnotesize{}1.1} & {\footnotesize{}0.1} & \tabularnewline
C-E & {\footnotesize{}1.85} & {\footnotesize{}2.0} & {\footnotesize{}3.4} & \textcolor{red}{\footnotesize{}5.0} & \textcolor{red}{\footnotesize{}5.0} & \textcolor{red}{\footnotesize{}5.0} & {\footnotesize{}3.8} & {\footnotesize{}2.1} & {\footnotesize{}0.6} & \tabularnewline
E-H & {\footnotesize{}1.7} & {\footnotesize{}2.0} & {\footnotesize{}3.4} & {\footnotesize{}4.8} & \textcolor{red}{\footnotesize{}5.0} & \textcolor{red}{\footnotesize{}5.0} & {\footnotesize{}3.8} & {\footnotesize{}2.3} & {\footnotesize{}0.9} & {\footnotesize{}0.1}\tabularnewline
C-F &  &  & {\footnotesize{}0.2} & {\footnotesize{}0.7} & {\footnotesize{}1.0} & {\footnotesize{}1.0} & {\footnotesize{}0.5} & {\footnotesize{}0.2} &  & \tabularnewline
F-H &  &  &  & {\footnotesize{}0.6} & {\footnotesize{}1.0} & {\footnotesize{}1.0} & {\footnotesize{}0.7} & {\footnotesize{}0.3} &  & \tabularnewline
\hline 
\end{tabular}
\par\end{centering}
\caption{\label{tab:Capacity-usage}Capacity usage. Capacity limitations marked
with red text color.}
\end{table}

The issue with small fractional flow volumes is not too problematic
when the number of traversed links (for each traffic relation) is
small. But as the number of links increase it becomes increasingly
problematic. We have unsuccessfully tried to figure out some additional
handling within the above variant of the multi-commodity formulation.
Instead we have chosen to redefine the layering structure of the network
so as to have one layer for each routing option of the traffic relations.
The input data for flow shares within the same time period ($\nu$)
is also extended so as to give the possible amount that can reach
every node along the specific route. The model formulation for this
is given in an accompanying document \textcite{aronsson_flows_2021}.

\end{document}

%% file: figures/bilevel_tikz.tex
%% Illustration of the bi-level model
\begin{tikzpicture}[auto]
  \node (upper) at (0,2) [rectangle, draw, text width=4cm, text centered] {Schedule upgrade projects};
  \node (lower) at (0,0) [rectangle, draw, text width=4cm, text centered] {Calculate traffic flow};
  \draw [->] (upper) to [out=0,in=0] node [rectangle, text width=3cm] {Capacity restrictions (TCR)} (lower);
  \draw [->] (lower) to [out=180,in=180] node {Adjusted traffic} (upper);
\end{tikzpicture}

%% file: figures/flowidea_tikz.tex
%% Illustration of basic flow idea
\begin{tikzpicture}[auto]
%% Time boxes and nodes
\foreach \x/\y in {0/1, 4/1, 0/-4, 4/-4}
	\draw (\x, \y) rectangle ++(4, 3);
\draw (2, 0) ellipse (2 and 0.3);
\draw (6, 0) ellipse (2 and 0.3);
%% Text markers
\draw (-0.5, 2.5) node[align=center] {Link\\$p$};
\draw (-0.5, -2.5) node[align=center] {Link\\$l$};
\draw (0, -4.5) node[anchor=west] {Time period $t$};
\draw (6, -4.5) node {$t+1$};
%% Train/speed lines
\foreach \x/\y/\dx in {0/4/1, 3/4/1, 4/4/1, 7/4/1,
					   0/-1/1.5, 2.5/-1/1.5, 4/-1/1.5, 6.5/-1/1.5}
	\draw (\x, \y) -- ++(\dx, -3);

%% Flow arcs
% Node inventory
\draw[->,very thick] (-0.2, 0) -- (0.2, 0);
\draw[->,very thick] (3.8, 0) -- (4.2, 0);
\draw[->,very thick] (7.8, 0) -- (8.2, 0);
\draw (8.2, 0) node[anchor=west] {$x^{Ni}_{n*}$};
% Transportation flow
\draw[->,very thick] (-0.5, 4.2) -- ++(1/3*4.1, -4.1);
\draw[->,very thick] (1.5, 4.2) -- ++(1/3*4.1, -4.1);
\draw[->,very thick] (3.5, 4.2) -- ++(1/3*4.1, -4.1);
\draw[->,very thick] (5.5, 4.2) -- ++(1/3*4.1, -4.1);
\draw[->,very thick] (0.9, -0.1) -- ++(1.5/3*4.1, -4.1);
\draw[->,very thick] (2.9, -0.1) -- ++(1.5/3*4.1, -4.1);
\draw[->,very thick] (4.9, -0.1) -- ++(1.5/3*4.1, -4.1);
\draw[->,very thick] (6.9, -0.1) -- ++(1.5/3*4.1, -4.1);
% with markers
\draw (-0.5, 4.2) node[anchor=south west] {$x^{Tr,f}_{p,t-1}$};
\draw (1.5, 4.2) node[anchor=south west] {$x^{Tr,d}_{pt}$};
\draw (3.5, 4.2) node[anchor=south west] {$x^{Tr,f}_{pt}$};
\draw (5.5, 4.2) node[anchor=south west] {$x^{Tr,d}_{p,t+1}$};
\draw (1.2, -0.6) node[anchor=west] {$x^{Tr,d}_{lt}$};
\draw (3.2, -0.6) node[anchor=west] {$x^{Tr,f}_{lt}$};
\draw (5.2, -0.6) node[anchor=west] {$x^{Tr,d}_{l,t+1}$};
\draw (7.2, -0.6) node[anchor=west] {$x^{Tr,f}_{l,t+1}$};

%% Flow share markers
\draw (0.05, 1.35) node[anchor=west, fill=white] {$\nu_p$};
\draw (1.5, 1.35) node[anchor=west, fill=white] {$1-\nu_p$};
\draw (4.05, 1.35) node[anchor=west, fill=white] {$\nu_p$};
\draw (5.5, 1.35) node[anchor=west, fill=white] {$1-\nu_p$};
\draw (1.25, -1.3) node[fill=white] {$1-\nu_l$};
\draw (3.25, -1.3) node[fill=white] {$\nu_l$};
\draw (5.25, -1.3) node[fill=white] {$1-\nu_l$};
\draw (7.25, -1.3) node[fill=white] {$\nu_l$};
% and cut-off time
\draw[dashed] (3, 4) -- ++(0, 1);
\draw[dashed] (7, 4) -- ++(0, 1);
\draw[dashed] (2.5, -1) -- ++(0, 2);
\draw[dashed] (6.5, -1) -- ++(0, 2);

\end{tikzpicture}

%% file: figures/ganttsmall_tikz.tex
%% Illustration two alternate renewal schedule solutions
\begin{tikzpicture}[auto]
%% Time tics
\draw (0.0, 0.0) node[anchor=south west] {Time period};
\foreach \i in {1,5,...,41} \draw (\i*0.25, 0) node {\i};
%% First solution
\draw (-1, -0.5) node {$|T|=35$};
\draw (-1, -0.8) node {P1};
\foreach \i/\x/\w in {1/14/2, 2/21/4, 3/29/2}
	\draw (\x*0.25, -1.0) rectangle ++(\w*0.25, 0.25) ++(-\w*0.125, -0.125) node [font=\tiny] {w\i};
\draw (-1, -1.1) node {P2};
\foreach \i/\x/\w in {1/1/2, 2/5/2, 3/9/2, 4/13/4, 5/21/4, 6/29/2, 7/33/2}
	\draw (\x*0.25, -1.3) rectangle ++(\w*0.25, 0.25) ++(-\w*0.125, -0.125) node [font=\tiny] {w\i};
\draw (-1, -1.4) node {P3};
\foreach \i/\x/\w in {1/2/1, 2/5/2, 3/9/8, 4/29/2, 5/33/1}
	\draw (\x*0.25, -1.6) rectangle ++(\w*0.25, 0.25) ++(-\w*0.125, -0.125) node [font=\tiny] {w\i};

%% Second solution
\draw (-1, -2.0) node {$|T|=40$};
\draw (-1, -2.3) node {P1};
\foreach \i/\x/\w in {1/17/2, 2/25/4, 3/34/2}
	\draw (\x*0.25, -2.5) rectangle ++(\w*0.25, 0.25) ++(-\w*0.125, -0.125) node [font=\tiny] {w\i};
\draw (-1, -2.6) node {P2};
\foreach \i/\x/\w in {1/1/2, 2/6/2, 3/11/2, 4/15/4, 5/25/4, 6/34/2, 7/38/2}
	\draw (\x*0.25, -2.8) rectangle ++(\w*0.25, 0.25) ++(-\w*0.125, -0.125) node [font=\tiny] {w\i};
\draw (-1, -2.9) node {P3};
\foreach \i/\x/\w in {1/2/1, 2/6/2, 3/11/8, 4/27/2, 5/35/1}
	\draw (\x*0.25, -3.1) rectangle ++(\w*0.25, 0.25) ++(-\w*0.125, -0.125) node [font=\tiny] {w\i};

\end{tikzpicture}

%% file: figures/gantt_resuse_tikz.tex
%% Illustration renewal schedules with minimal resource usage
\begin{tikzpicture}[auto]
%% Time tics
\draw (0.0, 0.0) node[anchor=south west] {Time period};
\foreach \i in {1,5,...,41} \draw (\i*0.25, 0) node {\i};
%% First solution
\draw (-1, -0.5) node {$|T|=35$};
\draw (-1, -0.8) node {P1};
\foreach \i/\x/\w in {1/2/2, 2/6/4, 3/17/2}
	\draw (\x*0.25, -1.0) rectangle ++(\w*0.25, 0.25) ++(-\w*0.125, -0.125) node [font=\tiny] {w\i};
\draw (-1, -1.1) node {P2};
\foreach \i/\x/\w in {1/1/2, 2/5/2, 3/9/2, 4/13/4, 5/21/4, 6/29/2, 7/33/2}
	\draw (\x*0.25, -1.3) rectangle ++(\w*0.25, 0.25) ++(-\w*0.125, -0.125) node [font=\tiny] {w\i};
\draw (-1, -1.4) node {P3};
\foreach \i/\x/\w in {1/1/1, 2/4/2, 3/10/8, 4/29/2, 5/33/1}
	\draw (\x*0.25, -1.6) rectangle ++(\w*0.25, 0.25) ++(-\w*0.125, -0.125) node [font=\tiny] {w\i};

%% Second solution
\draw (-1, -2.0) node {$|T|=40$};
\draw (-1, -2.3) node {P1};
\foreach \i/\x/\w in {1/22/2, 2/29/4, 3/38/2}
	\draw (\x*0.25, -2.5) rectangle ++(\w*0.25, 0.25) ++(-\w*0.125, -0.125) node [font=\tiny] {w\i};
\draw (-1, -2.6) node {P2};
\foreach \i/\x/\w in {1/1/2, 2/5/2, 3/9/2, 4/13/4, 5/25/4, 6/33/2, 7/37/2}
	\draw (\x*0.25, -2.8) rectangle ++(\w*0.25, 0.25) ++(-\w*0.125, -0.125) node [font=\tiny] {w\i};
\draw (-1, -2.9) node {P3};
\foreach \i/\x/\w in {1/4/1, 2/11/2, 3/17/8, 4/35/2, 5/39/1}
	\draw (\x*0.25, -3.1) rectangle ++(\w*0.25, 0.25) ++(-\w*0.125, -0.125) node [font=\tiny] {w\i};

%% Third solution
\draw (-1, -3.5) node {$|T|=45$};
\draw (-1, -3.8) node {P1};
\foreach \i/\x/\w in {1/2/2, 2/7/4, 3/15/2}
	\draw (\x*0.25, -4.0) rectangle ++(\w*0.25, 0.25) ++(-\w*0.125, -0.125) node [font=\tiny] {w\i};
\draw (-1, -4.1) node {P2};
\foreach \i/\x/\w in {1/1/2, 2/5/2, 3/11/2, 4/17/4, 5/29/4, 6/37/2, 7/43/2}
	\draw (\x*0.25, -4.3) rectangle ++(\w*0.25, 0.25) ++(-\w*0.125, -0.125) node [font=\tiny] {w\i};
\draw (-1, -4.4) node {P3};
\foreach \i/\x/\w in {1/5/1, 2/13/2, 3/21/8, 4/40/2, 5/44/1}
	\draw (\x*0.25, -4.6) rectangle ++(\w*0.25, 0.25) ++(-\w*0.125, -0.125) node [font=\tiny] {w\i};

\end{tikzpicture}

%% file: figures/gantt_both_tikz.tex
%% Illustration renewal schedules with minimal traffic disruption + resource usage
\begin{tikzpicture}[auto]
%% Time tics
\draw (0.0, 0.0) node[anchor=south west] {Time period};
\foreach \i in {1,5,...,41} \draw (\i*0.25, 0) node {\i};
%% First solution
\draw (-1, -0.5) node {$|T|=35$};
\draw (-1, -0.8) node {P1};
\foreach \i/\x/\w in {1/17/2, 2/21/4, 3/29/2}
	\draw (\x*0.25, -1.0) rectangle ++(\w*0.25, 0.25) ++(-\w*0.125, -0.125) node [font=\tiny] {w\i};
\draw (-1, -1.1) node {P2};
\foreach \i/\x/\w in {1/1/2, 2/5/2, 3/9/2, 4/13/4, 5/21/4, 6/29/2, 7/33/2}
	\draw (\x*0.25, -1.3) rectangle ++(\w*0.25, 0.25) ++(-\w*0.125, -0.125) node [font=\tiny] {w\i};
\draw (-1, -1.4) node {P3};
\foreach \i/\x/\w in {1/1/1, 2/5/2, 3/11/8, 4/29/2, 5/33/1}
	\draw (\x*0.25, -1.6) rectangle ++(\w*0.25, 0.25) ++(-\w*0.125, -0.125) node [font=\tiny] {w\i};

%% Second solution
\draw (-1, -2.0) node {$|T|=40$};
\draw (-1, -2.3) node {P1};
\foreach \i/\x/\w in {1/1/2, 2/5/4, 3/17/2}
	\draw (\x*0.25, -2.5) rectangle ++(\w*0.25, 0.25) ++(-\w*0.125, -0.125) node [font=\tiny] {w\i};
\draw (-1, -2.6) node {P2};
\foreach \i/\x/\w in {1/1/2, 2/5/2, 3/9/2, 4/13/4, 5/25/4, 6/33/2, 7/38/2}
	\draw (\x*0.25, -2.8) rectangle ++(\w*0.25, 0.25) ++(-\w*0.125, -0.125) node [font=\tiny] {w\i};
\draw (-1, -2.9) node {P3};
\foreach \i/\x/\w in {1/1/1, 2/11/2, 3/17/8, 4/35/2, 5/39/1}
	\draw (\x*0.25, -3.1) rectangle ++(\w*0.25, 0.25) ++(-\w*0.125, -0.125) node [font=\tiny] {w\i};

%% Third solution
\draw (-1, -3.5) node {$|T|=45$};
\draw (-1, -3.8) node {P1};
\foreach \i/\x/\w in {1/2/2, 2/5/4, 3/13/2}
	\draw (\x*0.25, -4.0) rectangle ++(\w*0.25, 0.25) ++(-\w*0.125, -0.125) node [font=\tiny] {w\i};
\draw (-1, -4.1) node {P2};
\foreach \i/\x/\w in {1/1/2, 2/9/2, 3/13/2, 4/17/4, 5/29/4, 6/37/2, 7/43/2}
	\draw (\x*0.25, -4.3) rectangle ++(\w*0.25, 0.25) ++(-\w*0.125, -0.125) node [font=\tiny] {w\i};
\draw (-1, -4.4) node {P3};
\foreach \i/\x/\w in {1/9/1, 2/15/2, 3/21/8, 4/39/2, 5/43/1}
	\draw (\x*0.25, -4.6) rectangle ++(\w*0.25, 0.25) ++(-\w*0.125, -0.125) node [font=\tiny] {w\i};

\end{tikzpicture}

%% file: figures/flowshare_tikz.tex
%% Illustration of flow share between current and next time period
\begin{tikzpicture}[auto]
%% Time boxes and nodes
\foreach \x/\y in {0/1, 4/1, 0/-4, 4/-4}
	\draw (\x, \y) rectangle ++(4, 3);
\draw (2, 0) ellipse (2 and 0.3);
\draw (6, 0) ellipse (2 and 0.3);
%% Text markers
%\draw (-0.5, 4.5) node {Links};
\draw (-0.5, 2.5) node[align=center] {Link\\$p$};
\draw (-0.5, -2.5) node[align=center] {Link\\$l$};
\draw (0, -4.5) node[anchor=west] {Time period $t$};
\draw (6, -4.5) node {$t+1$};
%% Train/speed lines
\foreach \x/\y/\dx in {0/4/1, 3/4/1, 4/4/1, 7/4/1,
					   0/-1/1.5, 2.5/-1/1.5, 4/-1/1.5, 6.5/-1/1.5}
	\draw (\x, \y) -- ++(\dx, -3);
%% Source/sink flow
\draw (0, 0.3) rectangle ++(4, 0.3);
\draw (4, 0.3) rectangle ++(4, 0.3);

%% Flow arcs
% Source/sink
\draw[->,very thick] (2, 0.45) -- ++(0, -0.3);
\draw[->,very thick] (6, 0.45) -- ++(0, -0.3);
\draw (8.2, 0.3) node[anchor=south west] {$x^{OD}_{*}$};
% Node inventory
%\draw (-0.2, 0) node[anchor=east] {$x^{Ni}_{n,t-1,h}$};
\draw[->,very thick] (-0.2, 0) -- (0.2, 0);
\draw[->,very thick] (3.8, 0) -- (4.2, 0);
\draw[->,very thick] (7.8, 0) -- (8.2, 0);
\draw (8.2, 0) node[anchor=west] {$x^{Ni}_{n*h}$};
% Transportation flow
\draw[->,very thick] (-0.5, 4.2) -- ++(1/3*4.1, -4.1);
\draw[->,very thick] (1.5, 4.2) -- ++(1/3*4.1, -4.1);
\draw[->,very thick] (3.5, 4.2) -- ++(1/3*4.1, -4.1);
\draw[->,very thick] (5.5, 4.2) -- ++(1/3*4.1, -4.1);
\draw[->,very thick] (0.9, -0.1) -- ++(1.5/3*4.1, -4.1);
\draw[->,very thick] (2.9, -0.1) -- ++(1.5/3*4.1, -4.1);
\draw[->,very thick] (4.9, -0.1) -- ++(1.5/3*4.1, -4.1);
\draw[->,very thick] (6.9, -0.1) -- ++(1.5/3*4.1, -4.1);
% with markers
\draw (-0.5, 4.2) node[anchor=south west] {$x^{Tr,n}_{p,t-1,h}$};
\draw (1.5, 4.2) node[anchor=south west] {$x^{Tr,d}_{pth}$};
\draw (3.5, 4.2) node[anchor=south west] {$x^{Tr,n}_{pth}$};
\draw (5.5, 4.2) node[anchor=south west] {$x^{Tr,d}_{p,t+1,h}$};
\draw (1.2, -0.6) node[anchor=west] {$x^{Tr,d}_{lth}$};
\draw (3.2, -0.6) node[anchor=west] {$x^{Tr,n}_{lth}$};
\draw (5.2, -0.6) node[anchor=west] {$x^{Tr,d}_{l,t+1,h}$};
\draw (7.2, -0.6) node[anchor=west] {$x^{Tr,n}_{l,t+1,h}$};

%% Flow share markers
\draw (0.05, 1.35) node[anchor=west, fill=white] {$\nu_p$};
\draw (1.5, 1.35) node[anchor=west, fill=white] {$1-\nu_p$};
\draw (4.05, 1.35) node[anchor=west, fill=white] {$\nu_p$};
\draw (5.5, 1.35) node[anchor=west, fill=white] {$1-\nu_p$};
\draw (1.25, -1.3) node[fill=white] {$1-\nu_l$};
\draw (3.25, -1.3) node[fill=white] {$\nu_l$};
\draw (5.25, -1.3) node[fill=white] {$1-\nu_l$};
\draw (7.25, -1.3) node[fill=white] {$\nu_l$};
% and cut-off time
\draw[dashed] (3, 4) -- ++(0, 1);
\draw[dashed] (7, 4) -- ++(0, 1);
\draw[dashed] (2.5, -1) -- ++(0, 2);
\draw[dashed] (6.5, -1) -- ++(0, 2);

\end{tikzpicture}